\documentclass{amsart}

\newtheorem{theorem}{Theorem}[section]
\newtheorem{lemma}[theorem]{Lemma}
\newtheorem{proposition}[theorem]{Proposition}
\newtheorem{corollary}[theorem]{Corollary}

\theoremstyle{definition}
\newtheorem{definition}[theorem]{Definition}
\newtheorem{example}[theorem]{Example}

\newtheorem{remark}[theorem]{Remark}
\newtheorem{construction}[theorem]{Construction}
\newtheorem*{acknowledgement}{Acknowledgement}

\theoremstyle{remark}

\usepackage{amscd,amssymb,graphics}
\usepackage{hyperref}

\newcommand\mylabel[1]{\label{#1}}

\newcommand{\cX}{{\underline{X}}}
\newcommand{\cU}{{\underline{U}}}
\newcommand{\cV}{{\underline{V}}}
\newcommand{\cW}{{\underline{W}}}
\newcommand{\cR}{{\underline{R}}}
\newcommand{\cS}{{\underline{S}}}
\newcommand{\cD}{{\underline{D}}}
\newcommand{\cC}{{\underline{C}}}
\newcommand{\cZ}{{\underline{Z}}}
\newcommand{\cT}{{\underline{T}}}
\newcommand{\cY}{{\underline{Y}}}

\newcommand{\NN}{\mathbb{N}}
\newcommand{\ZZ}{\mathbb{Z}}

\newcommand{\RR}{\mathbb{R}}

\newcommand{\PP}{\mathbb{P}}
\renewcommand{\AA}{\mathbb{A}}

\newcommand {\shA} {\mathcal{A}}
\newcommand {\shAut} {\mathcal{A}\!\text{\textit{ut}}}
\newcommand {\shB} {\mathcal{B}}
\newcommand {\shC} {\mathcal{C}}

\newcommand {\shExt} {\mathcal{E} \!\text{\textit{xt}}}

\newcommand {\shF} {\mathcal{F}}
\newcommand {\shG} {\mathcal{G}}
\newcommand {\shH} {\mathcal{H}}
\newcommand {\shHom} {\mathcal{H}\!\text{\textit{om}}}
\newcommand {\shI} {\mathcal{I}}
\newcommand {\shIsom} {\mathcal{I}\!\text{\textit{som}}}

\newcommand {\shK} {\mathcal{K}}
\newcommand {\shM} {\mathcal{M}}

\newcommand {\shN} {\mathcal{N}}

\newcommand {\shL} {\mathcal{L}}

\newcommand {\shT} {\mathcal{T}}

\newcommand {\shP} {\mathcal{P}}

\newcommand {\Aut} {\operatorname{Aut}}

\newcommand {\Br} {\operatorname{Br}}

\newcommand {\cH} {\check{H}}

\newcommand {\dd} {\hbox{\kern0.1em$/$\kern-0.7em $/$\kern0.1em}}

\newcommand {\Et} {{\operatorname{Et}}}

\newcommand {\Ext} {\operatorname{Ext}}

\newcommand {\gp} {{\operatorname{gp}}}

\newcommand {\Hom} {\operatorname{Hom}}

\newcommand {\id} {\operatorname{id}}

\newcommand {\interiour}  {\operatorname{int}}

\newcommand {\Isom} {\operatorname{Isom}}

\newcommand {\la} {\leftarrow}

\newcommand {\dirlim} {\varinjlim}

\newcommand {\LS} {\operatorname{LS}}

\newcommand {\lra} {\longrightarrow}

\renewcommand{\O} {\mathcal{O}}

\newcommand {\Pic} {\operatorname{Pic}}

\newcommand {\pr} {\operatorname{pr}}
\newcommand {\Proj} {\operatorname{Proj}}

\newcommand {\quadand} {\quad\text{and}\quad}

\newcommand {\ra} {\rightarrow}

\newcommand {\sep} {{\operatorname{sep}}}

\newcommand {\sqM} {\operatorname{\overline{\shM}}}

\newcommand {\sh} {{\operatorname{sh}}}
\newcommand {\Sing} {\operatorname{Sing}}

\newcommand {\Spec} {\operatorname{Spec}}

\def\mydate{\number\day\space\ifcase\month \or January\or February\or March\or
April\or May\or June\or July\or
August\or September\or October\or November\or December\fi \space\number\year}

\begin{document}

\title[Toroidal crossings and logarithmic structures]
{Toroidal crossings and logarithmic structures}

\author[Stefan Schroer]{Stefan Schr\"oer}
\address{Mathematisches Institut, Heinrich-Heine-Universit\"at,
40225 D\"usseldorf, Germany}
\curraddr{}
\email{schroeer@math.uni-duesseldorf.de}

\author[Bernd Siebert]{Bernd Siebert}
\address{Mathematisches Institut,
 Albert-Ludwigs-Universit\"at,
 Eckerstra\ss e 1,
 79104 Freiburg, Germany}
\email{bernd.siebert@math.uni-freiburg.de}

\subjclass{14D06, 14L32, 14J32}

\dedicatory{Second revision, 29 November 2004}

\begin{abstract}
We generalize Friedman's notion of d-semistability, which is a
necessary condition for spaces with normal crossings to admit
smoothings with regular total space. Our generalization deals with
spaces that locally look like the boundary divisor in Gorenstein
toroidal embeddings. In this situation, we replace d-semistability by
the existence of global log structures for a given gerbe of local log
structures. This leads to cohomological descriptions for the
obstructions, existence, and automorphisms of log structures.
We also apply  toroidal crossings to mirror
symmetry, by giving a duality construction involving toroidal crossing
varieties whose irreducible components are toric varieties.
This duality reproduces a
version of Batyrev's construction of mirror pairs for hypersurfaces
in toric varieties, but it applies to a larger class, including
degenerate abelian varieties.
\end{abstract}

\maketitle

\tableofcontents

\section*{Introduction}
Deligne and Mumford \cite{Deligne; Mumford 1969}  showed that any
curve  with normal crossing singularities deforms to a smooth
curve. This is no longer true for  higher dimensional spaces. Friedman
\cite{Friedman 1983} discovered that an obstruction for the existence
of smoothings with regular total space is an invertible sheaf on the
singular locus. He calls spaces with normal crossing
singularities \emph{d-semistable} if this sheaf is trivial.  So
d-semistability is a necessary condition for the existence of
smoothings with regular total space. This condition, however, is far
from being sufficient.

Nowadays, the notion of d-semistability is best understood via
\emph{log structures} and \emph{log spaces}. These  concepts are due to Fontaine and
Illusie, and were largely explored by K.\ Kato \cite{Kato 1989}.
They now play an important role in crystalline cohomology and
deformation theory, and have applications to Hodge theory, \'etale cohomology,
fundamental groups, and mirror symmetry.

Let $\cX$ be an algebraic space. A log space $X$ with underlying
algebraic space $\cX$ is a sheaf of monoids $\shM_X$ on $\cX$ related
to the structure sheaf $\O_\cX$ by certain axioms (for details see
Section \ref{spaces and structures}). It turns out that a space with
normal crossing singularities $\cX$ locally supports local log
structures that are unique up to isomorphism, and d-semistability is
equivalent with the existence of a global log structure, together
with the triviality of the normal sheaf. This has been exploited by
Kawamata and Namikawa  \cite{Kawamata; Namikawa 1994} for the
construction of Calabi--Yau manifolds, and by Steenbrink
\cite{Steenbrink 1995} for the construction of mixed Hodge
structures. The corresponding theoretical framework is due to Kato
\cite{Kato 1996}.

The first goal of this paper is to generalize the notion of
d-semistability to spaces that are locally isomorphic to boundary
divisors in toric varieties, which one might call \emph{toroidal
crossings}.  The theory of log structures suggests that such
generalization is possible, because spaces with normal crossing
singularities are just special instances of log smooth morphisms.
Furthermore, it became clear in various areas that spaces with normal crossing
singularities do not account for all degenerations that one wants to
study. Compare, for example, the work of Koll\'ar and Shepherd--Barron
\cite{Kollar; Shepherd-Barron 1988} on degenerate surfaces, and of
Alexeev \cite{Alexeev 2002} on degenerate abelian varieties.

Our main idea is to use the theory of \emph{nonabelian cohomology}, 
in particular the notion of \emph{gerbe}, to define d-semistability.
Roughly speaking, we define a \emph{log atlas} $\shG$ on $\cX$ to be a
gerbe of local log structures, that is, a collection of  locally
isomorphic log structures on certain neighborhoods. Now
d-semistability is nothing but the triviality of the gerbe class
$[\shG]$ in a suitable $H^2$-group, plus the triviality of the normal
bundle.

It turns out that the coefficient sheaf of the $H^2$-cohomology,
which is the \emph{band} in the terminology of nonabelian cohomology,
becomes abelian if we fix as additional datum the sheaf of monoids
$\sqM_X=\shM_X/\shM_X^\times$. Then the coefficient sheaf becomes the
abelian  sheaf $\shA_X=\shHom(\sqM_X,\O_\cX^\times)$, and this allows
us to calculate $H^2(\cX,\shA_X)$  via certain exact sequences. Using
such exact sequences, we deduce criteria for the existence of global
log structures. The theory works best if we make two additional
assumptions: The toric varieties that are local models should be
Gorenstein and regular in codimension two. Our main result is:
Each log atlas $\shG$ determines an invertible sheaf on
$\Sing(\cX)$ called the \emph{restricted conormal sheaf}, and its
extendibility to  $\cX$ is equivalent to the existence of a global log
structure, and its triviality is equivalent to d-semistability.
Olsson independently obtained
similar results in the case of normal crossing varieties \cite{Olsson 2003}.
He also showed that moduli of log structures yield
algebraic stacks \cite{Olsson 2003b}.

The second goal of this paper is to apply toroidal crossings to
mirror symmetry. Our starting point is the observation that from
$\overline{\shM}_X=\shM_X /\O_\cX^\times$ it is possible to construct
another degenerate variety $Y$, by gluing together the projective
toric varieties  $\Proj k[\overline\shM_{X,\overline x}^\vee]$ for
$x\in|\cX|$.  Furthermore, if $X$ itself consists of projective toric
varieties glued to each other along toric subvarieties, then there is
a monoid sheaf on $Y$ that at least locally is the sheaf $\sqM_Y$ of
a toroidal crossing log structure. This gives an involutive
correspondence between certain degenerate varieties endowed with such
sheaves of monoids.  Applied to hypersurfaces in projective toric
varieties it reproduces a degenerate version of Batyrev's mirror
construction, but it applies to many more cases, for example
degenerate abelian varieties.

Of course, mirror symmetry should do much more than what the naive
version presented here does. Our approach indicates that
one should try to understand mirror symmetry in terms of limiting
data of degenerations of varieties with trivial canonical bundle. By
\emph{limiting data} we mean  information about the degeneration
supported on the central fiber: Most importantly, the log structure induced by
the embedding into the total space, and certain cohomology classes on
the central fiber obtained by specialization. Mirror symmetry then  is
a symmetry acting on such limiting data. The
explanation of the mirror phenomenon would then be that it relates
limiting data of different degenerations. 
Mark Gross and the second author 
worked out a 
correspondence of true log spaces that involves also data
encoding the degeneration of a polarization,
see \cite{Gross; Siebert 2002} and \cite{Gross; Siebert 2003}.

\begin{acknowledgement}
We thank Mark Gross for valuable discussion about
the mirror construction in the last three sections, and
for finding  a mistake in our first version of Theorem \ref{canonical bijection}.
We also thank the referees for pointing out errors
and suggesting improvements.

Most of the work on this paper has been done while the first author
was a DFG research fellow at M.I.T., and the second author was visiting
the mathematical department of Jussieu as a DFG Heisenberg fellow.
We thank these institutions for hospitality.
Our project received financial support from the
DFG-Forschungsschwerpunkt \emph{Globale Methoden in der komplexen
Geometrie}. We thank the DFG for its support.
\end{acknowledgement}

\section{Algebraic spaces and logarithmic structures}
\label{spaces and structures}

In this section we recall some definitions regarding algebraic spaces
and logarithmic structures. For more details on algebraic spaces we
refer to the books of Knutson \cite{Knutson 1971} and Laumon and
Moret-Bailly \cite{Laumon; Moret-Bailly 2000}. For logarithmic
structures the reference is Kato's article \cite{Kato 1989}.
For typographical reasons, we use the following convention throughout:
Unadorned symbols $X,U,\ldots$ denote log spaces, whereas
underlined symbols $\cX,\cU,\ldots$ denote their underlying algebraic
spaces, as in \cite{Olsson 2003}. 

An \emph{algebraic space} $\cX$ is the quotient of a scheme $\cU$ by
an equivalence relation $\cR$, such that $\cR$ is a scheme, the
projections $\cR\ra \cU$ are \'etale, and the diagonal $\cR\ra
\cU\times \cU$ is a quasicompact monomorphism \cite{Knutson 1971}. 
Here quotient means quotient of sheaves on the site of rings endowed
with the \'etale topology. We prefer to work with algebraic spaces 
because operations like  gluing schemes yield algebraic spaces rather
than schemes (\cite{Artin 1970}, Theorem 6.1). Note that, over the
complex numbers, proper algebraic spaces correspond to compact
Moishezon spaces (\cite{Artin 1970}, Theorem 7.3).

A \emph{point} for $\cX$ is an equivalence class of morphisms
$\Spec(K)\ra \cX$, where $K$ is a field (\cite{Laumon; Moret-Bailly
2000}, Definition 5.2). The collection of all points is a topological
space $|\cX|$, whose open sets correspond to open subspaces
$\cU\subset \cX$. A morphism of algebraic spaces is called
\emph{surjective} if the induced map on the associated topological
spaces is surjective.

Let $\Et(\cX)$ be the \'etale site for $\cX$, whose objects are
the \'etale morphisms $\cU\ra \cX$, and whose covering families are
the surjections. 
A sheaf on $\cX$ is, by definition, a sheaf on
$\Et(\cX)$. Given a sheaf $\shF$ and a point $x\in|\cX|$, one defines
the stalk $\shF_{\bar{x}}=\dirlim\Gamma(\cU,\shF)$, where the
direct limit runs over all affine \'etale neighborhoods $\cU\ra \cX$
endowed with a point $u\in \cU$ such that
$\Spec\kappa(u)\ra \cX$ represents $x$.
Then $\shF\mapsto \shF_{\bar{x}}$ defines a fiber functor in the
sense of topos theory. According to \cite{SGA 4b}, Expos\'e
VIII, Theorem 3.5, a map between sheaves is bijective if and
only if for all points $x\in|\cX|$ the induced map between stalks are
bijective. Moreover, by \cite{SGA 4b}, Expos\'e VIII, Theorem 7.9,
the map
$x\mapsto\shF_{\bar{x}}$ is a homeomorphism between
$|\cX|$ and the space of topos-theoretical points for the topos
of sheaves on $\Et(\cX)$.

Let $\cX$ be an algebraic space. A \emph{log structure} on $\cX$ is
a sheaf of monoids $\shM_X$ on $\cX$ together with a homomorphism of monoids 
$\alpha_X:\shM_X\ra\O_{\cX}$ into the  multiplicative monoid $\O_{\cX}$, such that the induced map
$\alpha_X^{-1}(\O_{\cX}^\times)\ra \O_{\cX}^\times$ is
bijective \cite{Kato 1989}. 
A \emph{log space} $X$ is an algebraic space endowed with a log structure.
In other words, $X=(\cX,\shM_X,\alpha_X)$.


A \emph{chart} for a log space $X$ is an \'etale
neighborhood
$\cU\ra \cX$, together with a monoid $P$ and a homomorphism
$P\ra\Gamma(\cU,\O_\cX)$ so that the log space $U$ induced form the log space $X$ is isomorphic
to the log space associated to the constant \emph{prelog structure}
$P_{\cU}\ra\O_{\cU}$ (see \cite{Kato 1989}, Section 1 for details).
A log space is called \emph{fine} if it is covered by charts where the monoid
$P$ is fine, that is, finitely generated and integral.

Each log space
$X=(\cX,\shM_X,\alpha_X)$ comes along with a sheaf of
monoids
$$
\sqM_X=\shM_X/\shM_X^\times.
$$
Using the identifications 
$\shM_X^\times=\alpha_X^{-1}(\O_\cX^\times)=\O_\cX^\times$,
we usually write $\sqM_X=\shM_X/\O_\cX^\times$.
We call it
the \emph{ghost sheaf} of the log structure.
The stalks of the ghost sheaf are \emph{sharp} monoids, that is,
they have no units except the neutral element.
Ghost sheaves of fine log structures are not arbitrary.
Following \cite{SGA 4c}, Expos\'e IX, Definition 2.3,
we call a monoid sheaf $\shF$ \emph{constructible} if its stalks
are fine, and any affine \'etale neighborhood $\cU\ra \cX$
admits a decomposition into finitely many constructible locally closed
subschemes $\cU_i$ such that the restrictions $\shF_{\cU_i}$ are locally
constant.

\begin{proposition}
\mylabel{constructible}
If $X$ is a fine log space, then its ghost sheaf $\sqM_X$ is a
constructible monoid sheaf.
\end{proposition}

\proof
This is a local problem by
\cite{SGA 4c}, Expos\'e IX, Proposition 2.8.
Hence we easily reduce to the case that
$X=\Spec(\ZZ[P])$ for some fine monoid $P=\sum_{i=1}^r\ZZ p_i$.
Each subset $J\subset\left\{1,\ldots,r\right\}$ yields a  ring $R_J=S^{-1}\ZZ[P]/I$,
where $S\subset \ZZ[P]$ is the multiplicative subset generated
by all $p_i$ with $i\not\in J$, and $I\subset S^{-1}\ZZ[P]$
is the ideal generated by all $p_i$ with $i\in J$.
We obtain locally closed subsets $\cX_J=\Spec(R_J)$.
Note that $\cX_J\subset \cX$ is the set of points $x\in \cX$
where the sections $p_i$ are invertible if $i\not\in J$,
and vanish if $i\in J$.
It follows that we have a  disjoint decomposition 
$\cX=\bigcup_J \cX_J$.

To see that $\sqM_X$ is constant along $\cX_J$, fix a point 
$x\in \cX_J$.
Then the germ $\sqM_{X,\bar x}$ equals the sharp monoid
$(P+\sum_{i\not\in J}\ZZ p_i)/G$, where $G\subset P+\sum_{i\not\in J}\ZZ p_i$
is the subgroup of invertible elements. This does not depend on the
point $x$, hence the assertion.
\qed

\medskip
Given an algebraic space $\cX$ and  two points $x,y\in|\cX|$ with $y\in\overline{\left\{x\right\}}$,
one has a specialization map $\shF_{\bar{y}}\ra\shF_{\bar{x}}$
(some authors call it a cospecialization map).
We say that $\shF$
has \emph{surjective specialization maps} if these maps are surjective
for all pairs $x,y\in|\cX|$ with $y\in\overline{\left\{x\right\}}$.
Ghost sheaves are typical examples:

\begin{proposition}
\mylabel{surjective specialization}
Let $X$ be a log space. If each point $x\in|\cX|$ admits a chart,
then the ghost sheaf $\sqM_X$ has surjective specialization maps.
\end{proposition}

\proof
This is a local problem, and we may assume that
$\cX=\Spec(A)$ is the spectrum of a henselian local ring with separably
closed residue field.
Choose a monoid $P$ and a map $f:P\ra\Gamma(\cX,\O_\cX)$ so that
$X$ is the associated log space.
The cocartesian diagram
$$
\begin{CD}
f^{-1}(\O_\cX^\times) @>>> P_\cX\\
@VVV @VVV\\
\O^\times_\cX @>>> \shM_X
\end{CD}
$$
yields the monoid sheaf $\shM_X$. As a consequence,
the composite map $P_\cX\ra\sqM_X$ is surjective, and the ghost sheaf
$\sqM_X$ has surjective specialization maps.
\qed

\section{Logarithmic atlases}

Let $\cX$ be an algebraic space. A natural question to ask is:
What is the set of \emph{all} log spaces $X$ with underlying algebraic
space $\cX$, up to isomorphism?
This is a reasonable moduli problem, as
Olsson \cite{Olsson 2003b} proved that the fibered category of
fine log structures on $\cX$-schemes $\cU$ is  an algebraic
stack.
Here we seek a cohomological approach to classify log structures.
This classification problem, however, is nonabelian in nature. To overcome this, we shall fix the ghost sheaf $\sqM_X$ such that
the problem becomes abelian. This leads to the desired cohomological
descriptions for obstructions, existence, and automorphisms of log
structures.

Given a log space $X$ with underlying algebraic space $\cX$,
we denote by $\shAut_{X/\cX}$ the sheaf of log space automorphisms
$X\ra X$ inducing the
identity on the underlying algebraic space $\cX$. Such automorphism
correspond to bijections $\phi:\shM_X\ra\shM_X$  compatible
with $\alpha_X:\shM_X\ra\O_\cX$. They necessarily fix
the subsheaf $\O_\cX^\times\subset\shM_X$ pointwise, and induce a bijection
$\overline{\phi}:\sqM_X\ra\sqM_X$. Let
$\shAut'_{X/\cX}\subset\shAut_{X/\cX}$ be the subsheaf of automorphisms
inducing the identity on the ghost sheaf $\sqM_X$.
We want to compare $\shAut'_{X/\cX}$ to the abelian sheaf
$$
\shA_X= \shHom(\sqM_X,\O_\cX^\times)=\shHom(\sqM_X^\gp,\O_\cX^\times).
$$
There is a canonical inclusion
$\shA_X\subset\shAut(\shM_X)$ sending a map
$h:\sqM_X\ra\O_\cX^\times$ to
$$
\shM_X\lra\shM_X,\quad s\longmapsto s+ h(\overline{s}),
$$
where $\overline{s}\in\Gamma(\cU,\sqM_X)$ denotes the image of
$s\in\Gamma(\cU,\shM_X)$, and $\cU\ra \cX$ is any affine \'etale neighborhood.

\begin{proposition}
Suppose $\cX$ is reduced. Then the inclusion
$\shA_X\subset\shAut(\shM_X)$ factors over the inclusion
$\shAut'_{X/\cX}\subset\shAut(\shM_X)$.
\end{proposition}

\proof
With the preceding notation, we have to check
that equality
$$\alpha_X(s)\cdot\alpha_X(h(\overline{s}))=\alpha_X(s)
$$
holds inside
$\Gamma(\cU,\O_\cX)$.
This is obvious if $\alpha_X(s)=0$. If $\alpha_X(s)$
is invertible  then $\overline{s}=0$, and equality holds as well.
Let $\eta_i\in\cU$, $i\in I$ be the generic points.
Since $\cX$ is reduced, there are open neighborhoods $\eta_i\in \cU_i$
so that $\alpha(s)_{\cU_i}$ is either zero or invertible.
We infer that the desired equality holds on $\bigcup_{i\in I}\cU_i$.
Using again that $\cX$ is reduced, we see that
$\alpha_X(s)\cdot\alpha_X(h(\overline{s}))=\alpha_X(s)$ holds on $\cU$.
\qed

\begin{proposition}
\mylabel{automorphism sheaf}
Suppose $\cX$ is reduced and  $\sqM_X$ has integral stalks. Then the inclusion
$\shA_X\subset\shAut'_{X/\cX}$ is bijective.
\end{proposition}

\proof
Fix a point $x\in |\cX|$. We have to show that the inclusion
$\shA_{X,\bar{x}}\subset\shAut'_{X/\cX,\bar{x}}$ is bijective.
Let $\cU\ra \cX$ be an
\'etale neighborhood of
$x$ and $U$ the induced log space, and $\phi:\shM_U\ra\shM_U$ a bijection compatible
with $\alpha_U$ and inducing
the identity on
$\sqM_U$. We now construct a homomorphism $h:\sqM_U\ra\O_U^\times$ as
follows:

Let $\bar{s}\in\Gamma(\cV,\sqM_X)$ be a local section on an \'etale
neighborhood
$\cV\ra \cU$. Choose a refinement $\cW\ra \cV$ so that $\bar{s}_{\cW}$ comes from
a section $s\in\Gamma(\cW,\shM_X)$. Then the equation
$\phi(s)=s +\alpha^{-1}(t)$ defines a section
$t\in\Gamma(\cW,\O_\cX^\times)$. Since $\sqM_X$ has integral stalks,
so has
$\shM_X$, and we infer from the defining equation that $t$
depends only on
$\overline{s}_{\cW}$, and not on the choice of $s$.
Consequently 
$\pr_0(t)=\pr_1^*(t)$ on $\cW\times_{\cV} \cW$, and
$t$ descends to a section $h(\overline{s})\in\Gamma(\cV,\O_\cX^\times)$.
Furthermore, this
section depends only on
$\overline{s}$, and not on the choice of the refinement
$\cW\ra\cV$.

It follows from the defining equation $\phi(s)=s + \alpha^{-1}(t)$
that $h(\overline{s})$ yields a monoid homomorphism
$\Gamma(\cV,\sqM_X)\ra\Gamma(\cV,\O_\cX^\times)$. Clearly
$h(\overline{s})$ is compatible with restrictions. Hence we have defined a
sheaf homomorphism
$h:\sqM_U\ra\O_U^\times$, with $\phi(s)=s+ h(\overline{s})$ for any local
section
$s\in\Gamma(\cV,\shM_X)$. In other words, the germ
$h_{\bar{x}}\in\shA_{X,\bar{x}}$ corresponds to the germ
$\phi_{\bar{x}}\in\shAut'_{X/\cX,\bar{x}}$
under the canonical inclusion.
\qed

\medskip
Now let $\cX$ be an algebraic space, and fix as additional datum
a sheaf of integral sharp monoids $\sqM_{\cX}$.
Let $\LS(\cX)$ be the category of pairs $(U,\varphi)$,
where $U=(\cU,\shM_U,\alpha_U)$
is a log space, whose underlying algebraic space $\cU$ is an \'etale neighborhood $\cU\ra\cX$,
and
$$
\varphi:\sqM_U=\shM_U/\O_U^\times\lra\sqM_{\cX}|_{\cU}=\sqM_{\cU}
$$
is an isomorphism. We call $\varphi$ a \emph{framing} for
the log space $U$ with respect to $\sqM_{\cX}$. The functor
$$
\LS(\cX)\lra\Et(\cX),\quad (U,\varphi)\longmapsto\cU
$$
yields a
fibered category.
The fiber $\LS(\cX)_{\cU}$ over an \'etale neighborhood $\cU$
is equivalent to the category of log structures on $\cU$
whose ghost sheaf is identified with $\sqM_{\cU}=\sqM_\cX|_{\cU}$. By abuse of notation,
we usually write $U$ instead of $(U,\varphi)$ for
the objects in $\LS(\cX)$. An inverse image for an $\cX$-morphism
of \'etale neighborhoods $g:\cU\ra \cV$ is given by restriction.
This also extends from the small \'etale site $\Et(\cX)$ to the big
\'etale site, where the preimage is given by
the log structure associated to the prelog structure
$g^{-1}(\shM_{V})\ra\O_{\cU}$.
Obviously, our fibered category is a \emph{stack} in Giraud's
sense \cite{Giraud 1971}, Chapter II, Definition 1.2,
that is, all descent data are effective.

Now recall that a substack
$\shG\subset \LS(\cX)$ over $\Et(\cX)$ is a
\emph{subgerbe} if, for each \'etale neighborhood $\cU\ra\cX$,
the following
axioms hold (see \cite{Giraud 1971}, Chapter III, Definition 2.1.3):
\renewcommand{\labelenumi}{(\roman{enumi})}
\begin{enumerate}
\item The objects in $\shG_{\cU}$ are locally isomorphic.
\item The morphisms in $\shG_{\cU}$ are isomorphisms.
\item There is an \'etale covering $\cV\ra \cU$ with
$\shG_\cV$ nonempty.
\end{enumerate}
A gerbe with $\shG_{\cX}\neq\emptyset$ is called \emph{neutral}.
This means that it is possible to glue the local
log structures $V\in \shG$, which exists by axiom
(iii), in at least one way to obtain a global log structure
$X\in\shG$. Note that,
with respect to inclusion, 
each subgerbe is contained in a maximal subgerbe, and
we may restrict our attention to maximal subgerbes. 
The following definition is fundamental for the rest of this paper:

\begin{definition}
\mylabel{log atlas}
Let $\cX$ be an algebraic space endowed with a sheaf $\sqM_{\cX}$
of integral sharp monoids.
A \emph{log atlas} for $\cX$ with respect to $\sqM_{\cX}$
is a maximal subgerbe $\shG\subset \LS(\cX)$ over $\Et(\cX)$.
\end{definition}

The idea is that a log atlas $\shG$ tells us how local log structures on
$\cX$ should be locally around each point $x\in |\cX|$, up to isomorphism.
It does not, however, single out preferred local log structures.
Neither does it inform us how to glue these local log structures.
Given a log atlas, the problem is to decide whether or not
it admits a global log structure.
Note that Kawamata and Namikawa \cite{Kawamata; Namikawa 1994}
used the word log atlas in a very different way, namely to denote
global log structures.

Given an object $(U,\varphi)\in\shG$, we obtain a
homomorphism
$$
\shA_{\cU}=\shHom(\sqM_{\cU},\O_U^\times) \stackrel{\varphi^*}{\lra}
\shHom(\sqM_U,\O_U^\times)\lra\shAut'_{U/\cU},
$$
which is bijective by Proposition \ref{automorphism sheaf}.
In the language of nonabelian cohomology, the abelian sheaf
$$
\shA_{\cX}=\shHom(\sqM_{\cX},\O_{\cX}^\times)
$$
binds the gerbe $\shG$, and $\shG$ becomes an $\shA_{\cX}$-gerbe
(\cite{Giraud 1971}, Chapter IV, Definition 2.2.2). In turn, we
obtain a gerbe class $[\shG]\in H^2(\cX,\shA_{\cX})$.
The theory of nonabelian cohomology immediately gives
the following:

\begin{theorem}
\mylabel{cohomological}
Let $\shG$ be a log atlas on an algebraic space $\cX$
with respect to a sheaf $\sqM_{\cX}$ of integral sharp monoids.
Then there is a global log structure $X\in\shG$
if and only if the gerbe class $[\shG]\in H^2(\cX,\shA_{\cX})$ vanishes.
In this case, the set
of isomorphism classes of $X\in\shG$
is a torsor for $H^1(\cX,\shA_{\cX})$.
Moreover, for each global log structure
$X\in\shG$, the group
of log space automorphisms inducing the identity on the underlying algebraic space $\cX$ and on
the sheaf
$\sqM_{\cX}$ is
$H^0(\cX,\shA_{\cX})$.
\end{theorem}

\proof
The first statement is \cite{Giraud 1971}, Chapter IV, Theorem 3.4.2.
The second statement follows from 
\cite{Giraud 1971}, Chapter III, Theorem 2.5.1. The last statement is
nothing but Proposition \ref{automorphism sheaf}.
\qed

\medskip
The preceding result is almost a tautology if we use
the \emph{geometric definition} for the universal
$\partial$-functor $H^n(\cX,\shF)$, where $\shF$ is an abelian sheaf
and $0\leq n\leq 2$.
In this definition, 
$H^1(\cX,\shF)$ is the set of isomorphism classes of $\shF$-torsors,
and $H^2(\cX,\shF)$ is the set of equivalence classes of
$\shF$-gerbes. Given a short exact sequence
$0\ra\shF'\ra\shF\ra\shF''\ra 0$, the coboundary operator maps a
section for $\shF''$ to the $\shF'$-torsor of its preimage in
$\shF$, and an
$\shF''$-torsor to the $\shF'$-gerbe of its $\shF$-liftings.

With these definitions, the cohomology class $[\shG]\in
H^2(\cX,\shA_{\cX})$ of a log atlas
$\shG$ is the equivalence class of the underlying $\shA_{\cX}$-gerbe,
and the difference between two isomorphism classes of global log
spaces
$X,X'\in\shG$ is the isomorphism class
of the $\shA_{\cX}$-torsor
$\shIsom(X', X)$.
The situation becomes more illuminating if we use other
descriptions for cohomology. We discuss this in the next section.

\section{Cohomology and hypercoverings}
\mylabel{cech cohomology}

Let us recall the \emph{combinatorial definition}
for cohomology in degrees $\leq 2$.
Let $\shF$ be an abelian sheaf on an algebraic space $\cX$. Then one may describe  $H^n(\cX,\shF)$ for $0\leq n\leq
2$ as follows.

Suppose we have an \'etale covering
$\cU\ra \cX$ and an
\'etale covering $\cV\ra \cU\times_{\cX} \cU$.
Let
$p_0:\cV\ra \cU$ be the composition $v\mapsto(u_0,u_1)\mapsto u_0$, and
$p_1:\cV\ra \cU$ be the other composition $v\mapsto(u_0,u_1)\mapsto u_1$.
Define
$$
(\cV/\cU)_2\subset \cV\times_{\cX}\cV\times_{\cX}\times \cV$$
to be the subspace
satisfying the simplicial identities
$p_i\circ \pr_j=p_{j-1}\circ \pr_i$,
$i<j$, and let $p_j:(\cV/\cU)_2\ra V$ be the maps induced by
the projections
$\pr_j$. Here $\pr_j$ denotes the projections 
$(v_0,v_1,v_2)\mapsto v_j$. Inductively, we define for each
$n\geq 2$ subspaces
$(\cV/\cU)_{n+1}\subset\prod_{i=0}^n(\cV/\cU)_n$ and projections
$p_j:(\cV/\cU)_{n+1}\ra (\cV/\cU)_n$ as above.
This gives a semisimplicial \'etale covering $(\cV/\cU)_\bullet$ of $\cX$,
where $(\cV/\cU)_1=\cV$ and $(\cV/\cU)_0=\cU$. In fact, $(\cV/\cU)_\bullet$ is the
coskeleton for the truncated semisimplicial covering
$\cV\rightrightarrows \cU$ (for more on this, see
\cite{Duskin 1975}, Section 0.7).

\begin{remark}
\mylabel{hypercovering etale}
The maps $p_j:(\cV/\cU)_{n+1}\ra (\cV/\cU)_{n}$ are indeed \'etale.
To see this, note first that the composite maps
$p_ip_j:(\cV/\cU)_{n+1}\ra (\cV/\cU)_{n-1}$ are \'etale, because $(\cV/\cU)_{n+1}$ is
defined as a fiber product with respect to \'etale maps.
By induction on $n$, the maps $p_i:(\cV/\cU)_{n}\ra (\cV/\cU)_{n-1}$ are \'etale, and it then
follows from \cite{EGA IVd}, Corollary 17.3.5 that $p_j:(\cV/\cU)_{n+1}\ra
(\cV/\cU)_{n}$ are \'etale as well.
\end{remark}

Now let $\shF$ be any abelian sheaf on $\cX$.
In accordance with the applications we have in mind,
we shall write the group law
multiplicatively. The sheaf $\shF$ yields a cochain complex of abelian groups
$
C^n(\cV/\cU,\shF)=\Gamma((\cV/\cU)_n,\shF)
$
with the usual differential
$d=\prod {p_i^*}^{(-1)^i}$. Let
$H^n(\cV/\cU,\shF)$ be the corresponding
cohomology group. Given other \'etale coverings
${\cV}{}'\rightrightarrows \cU{}'$ refining
the given \'etale coverings $\cV\rightrightarrows \cU$ , we obtain an
induced map
$H^n(\cV/\cU,\shF)\ra H^n(\cV{}'/\cU{}',\shF)$. Now let us \emph{define}
$$
H^n(\cX,\shF)=\dirlim H^n(\cV/\cU,\shF),\quad 0\leq n\leq 2.
$$
This is a $\partial$-functor: Given a short exact sequence
$0\ra\shF'\ra\shF\ra\shF''\ra 0$ and a 0-cocycle $f\in
Z^0(\cV/\cU,\shF'')$, we refine $\cU$, choose an
$\shF$-valued 0-cochain
lift $\tilde{f}$ of $f$, and define
$\partial(f)=p_0^*(\tilde{f})/p_1^*(\tilde{f})$. Similarly, given a
1-cocycle $g\in Z^1(\cV/\cU,\shF'')$, we pass to a refinement of $\cV$,
choose an
$\shF$-valued 1-cochain
lift $\tilde{g}$ of $g$, and define
$\partial(g)=p_0^*(\tilde{g})p_2^*(\tilde{g})/p_1^*(\tilde{g})$.
It is not difficult to see that this
$\partial$-functor vanishes on injective
sheaves, hence is universal by
\cite{Grothendieck 1957}, Proposition 2.2.1. Therefore, the
geometric and combinatorial definitions for
$H^n(\cX,\shF)$, $0\leq n\leq 2$ are canonically
isomorphic as $\partial$-functors.

The canonical isomorphism between geometric and combinatorial
definition takes the following explicit form:
Suppose we have an $\shF$-torsor $\shT$.
Choose an \'etale covering $\cU\ra \cX$ so that there is a section
$s\in\Gamma(\cU,\shT)$ and set $\cV=\cU\times_{\cX} \cU$. Then $p_0^*(s)=p_1^*(s)\cdot f$
defines a cocycle
$f\in Z^1(\cV/\cU,\shF)$. To see that $\shT\mapsto f$ yields the canonical
isomorphism, it suffices to check that the induced map is well-defined,
additive, and commutes 
with the coboundary
$\partial:H^0\ra H^1$, which is straightforward.

Now suppose $\shG$ is an $\shF$-gerbe. Choose an \'etale covering
$\cU\ra \cX$ admitting an object $T\in\shG_{\cU}$, and an \'etale covering
$\cV\ra \cU\times_{\cX} \cU$ admitting
an isomorphism $\phi:p_1^*(T)\ra p_0^*(T)$. Then the equation
$$
g\cdot p_1^*(\phi)=p_0^*(\phi)p_2^*(\phi)\in
\Isom(p_1^*p_1^*T, p_0^*p_0^*T)
$$
defines a cocycle $g\in Z^2(\cV/\cU,\shF)$. Note that this equation involves
the simplicial identities $p_j^*p_i^*(T)\simeq p_i^*p_{j-1}^*(T)$, $i<j$. To
see that
$\shG\mapsto g$ yields the canonical isomorphism, it suffices to check
that the induced map is well-defined,
additive, and commutes with
$\partial:H^1\ra H^2$, which is again straightforward.

The action of $H^1(\cX,\shF)$ on the set of isomorphism classes of
$\shG_{\cX}$ is as follows: Given a $\shF$-torsor $\shT$, choose a cocycle
$f\in Z^1(\cV/\cU,\shF)$ as above and a global object $T\in\shG_{\cX}$.
We have a canonical bijection $\phi:p_1^*(T_{\cU})\ra p_0^*(T_{\cU})$ on
$\cV=\cU\times_{\cX} \cU$. Then the isomorphism
$\phi\circ f:p_1^*(T_{\cU})\ra p_0^*(T_{\cU})$ is another descent datum,
that is,
\begin{equation}
\label{descent datum}
p_1^*(\phi)\circ p_1^*(f)=
p_0^*(\phi)\circ p_0^*(f)\circ p_2^*(\phi)\circ p_2^*(f)
\end{equation}
holds as isomorphisms  on $(\cV/\cU)_2$, with  suitable identifications
coming from the simplicial identities.
Indeed, we have $p_0^*(f)\circ p_2^*(\phi)=p_2^*(\phi)\circ p_0^*(f)$,
because $\shF$ is abelian, and (\ref{descent datum}) follows
from the cocycle condition for $\phi$ and $f$.
Summing up, the descend datum $\phi\circ f$
defines another global object
$T'\in\shG_{\cX}$, together with a bijection $\shT\ra\shIsom(T,T')$.

\begin{remark}
\mylabel{Cech or sheaf}
Note that we obtain \v{C}ech cohomology groups $\cH^n(\cX,\shF)$
if we use $\cV=\cU\times_{\cX} {\cU}$ instead of \'etale coverings
$\cV\ra \cU\times_{\cX} \cU$. In general
\v{C}ech cohomology groups do not form a $\partial$-functor
on the category of sheaves and differ from 
true cohomology groups. Note, however,
that the canonical map $\cH^n(\cX,\shF)\ra H^n(\cX,\shF)$ is bijective
for all $n\geq 0$
provided that $\cX$ admits an ample invertible sheaf \cite{Artin 1971}.
Furthermore, $\cH^2(\cX,\shF)\ra H^2(\cX,\shF)$ is bijective if each
pair of points admits an affine open neighborhood \cite{Schroeer 2002}.
\end{remark}

\section{The sheaf of automorphisms}
\mylabel{sheaf of automorphisms}

Let $\cX$ be an algebraic space, endowed with a sheaf of
integral sharp monoids $\sqM_\cX$. As before, we set
$\shA_\cX=\shHom(\sqM_\cX,\O_\cX^\times)$.
The goal now is to compute the cohomology groups $H^1(\cX,\shA_\cX)$ and
$H^2(\cX,\shA_\cX)$ in some interesting special cases. To this end we shall relate the sheaf
$\shA_\cX$ to other sheaves via exact
sequences. This relies on the following construction.

Suppose that our algebraic space $\cX$ is a noetherian, reduced, and
satifies the following condition:
For all points $x\in|\cX|$, the integral components of
$\Spec(\O_{\cX,\bar{x}})$ are
normal. This condition holds for the
spaces we have in mind for applications, namely
boundary divisors in toroidal embeddings.
The referees  pointed out that such a condition
is indeed indispensable.
The assumption implies that the normalization $f:\cS\ra \cX$ is a finite map.
Moreover, $f$ is an isomorphism near each point
$x\in|\cX|$ where $\O_{\cX,\bar{x}}$
is unibranch, because then $\O_{\cX,\bar{x}}$ is integral by
\cite{EGA IVd}, Corollary 18.6.13.

Let $\shI\subset\O_\cX$ be the conductor ideal for $f$, that is
the annihilator
ideal of $f_*(\O_\cS)/\O_\cX$, or equivalently the largest coherent $\O_\cX$-ideal
that is also an $\O_\cS$-ideal.
The closed subspaces $\cD\subset \cX$ and $f^{-1}(\cD)\subset \cS$
defined by
the conductor ideal are the branch space and the ramification space
for the finite morphism $f:\cS\ra \cX$, respectively.
We call $\cD\subset \cX$ the \emph{subspace of nonnormality}.
The cartesian diagram
$$
\begin{CD}
f^{-1}(\cD) @>>> \cS \\
@VVV @VV fV\\
\cD @>>> \cX
\end{CD}
$$
yields sequences of coherent $\O_\cX$-modules
\begin{equation}
\mylabel{additive}
0\lra\O_\cX\lra \O_\cS\oplus\O_\cD\lra\O_{f^{-1}(\cD)}\lra 0.
\end{equation}
Here the map on the left is the diagonal map $t\mapsto (t,\overline{t})$,
and the map on the right is the difference map $(t,\overline{s})\mapsto \overline{t}-\overline{s}$.
Similarly, we have a sequence of  abelian sheaves on $\cX$
\begin{equation}
\mylabel{multiplicative}
1\lra\O_\cX^\times\lra \O_\cS^\times\times\O_\cD^\times\lra\O_{f^{-1}(\cD)}^\times\lra 1.
\end{equation}
For the sake of simplicity we have supressed $f_*$ from notation.

\begin{proposition}
\mylabel{cocartesian}
The preceding sequences (\ref{additive}) and (\ref{multiplicative})
are exact.
\end{proposition}

\proof
This is a local problem, so we may assume that our algebraic spaces $\cX=\Spec(A)$ and
$\cS=\Spec(B)$ are affine. Let $I\subset A$ be the conductor ideal.
We treat the additive sequence (\ref{additive}) first.
It is easy to see that this sequence is a complex, and exact at $\O_\cX$ and
$\O_{f^{-1}(\cD)}$. To see that the complex is exact in the middle, suppose we have $(t,\overline{s})\in B\oplus A/I$ with
$\overline{t}=\overline{s}$. Subtracting the image of $s\in A$, we may
assume that $\overline{s}=0$. It then follows $t\in I\subset A$, so
$(t,0)$ lies in the image of the diagonal map $A\ra B\times A/I$.

It remains to treat the multiplicative sequence (\ref{multiplicative}).
Again it is immediate that this sequence is a complex that is exact at the outer terms.
To see that the complex is exact in the middle, suppose we have a pair $(t,\overline{s})\in B^\times\times (A/I)^\times$
with $\overline{t}/\overline{s}=1$. Then $\overline{s}=\overline{t}$,
and we just saw that this implies $t\in A$.
Repeating this argument with $(1/t,1/\overline{s})$, we see
that $1/t\in A$, hence $t\in A^\times$.
\qed

\medskip
Next, consider the constructible sheaf $f_*(\ZZ_\cS)$ on $\cX$.
Each stalk $f_*(\ZZ_\cS)_{\bar{x}}$ is a free $\ZZ$-modules whose rank
is the number of irreducible components in $\Spec(\O_{\cX,\bar{x}})$.
Let $\rho\in\Gamma(\cX,f_*(\ZZ_\cS))$ be the diagonal section defined by
$\rho_{\bar{x}}=(1,\ldots,1)$, which corresponds to
$1\in\Gamma(\cS,\ZZ)$. We have an evaluation map
$$
\rho^*:\shHom(f_*(\ZZ_\cS),\O_\cX^\times)\lra\O_\cX^\times, \quad s\longmapsto s(\rho)
$$
and a sequence of abelian sheaves
\begin{equation}
\label{sequence}
1\lra
\shHom(f_*(\ZZ_\cS),\O_\cX^\times)\stackrel{\rho^*}{\lra}\O_\cX^\times\lra
i_*(\O_\cD^\times)\lra 1,
\end{equation}
where $i:\cD\ra \cX$ denotes the closed embedding of the space of
nonnormality.

\begin{proposition}
\mylabel{exactness}
The preceding sequence (\ref{sequence}) is exact.
\end{proposition}

\proof
For simplicity we set $\shB=\shHom(f_*(\ZZ_\cS),\O_\cX^\times)$.
The short exact sequence
$0\ra\ZZ_\cX\stackrel{\rho}{\ra} f_*(\ZZ_\cS)\ra\shF\ra 0$ defines an abelian
sheaf
$\shF$ with $i^{-1}(\shF)=0$ for some dense open embedding $i:\cU\ra \cX$.
Applying $\shHom(.,\O_\cX^\times)$, we obtain an exact sequence
$$
1\lra\shHom(\shF,\O_\cX^\times)
\lra\shB_\cX\stackrel{\rho*}{\lra}\O_\cX^\times.
$$
There is an inclusion
$\shHom(\shF,\O_\cX^\times)\subset\shHom(\shF,i_*i^{-1}(\O_\cX^\times))$ because
$\cX$ has no embedded components. Moreover,
$\shHom(\shF,i_*i^{-1}(\O_\cX^\times))=
i_*\shHom(i^{-1}(\shF),i^{-1}(\O_\cX^\times))$
by \cite{SGA 2}, Expos\'e 1, Corollary 1.5.
The latter sheaf vanishes because $i^{-1}(\shF)=0$, and we conclude that
$\rho^*:\shB_\cX\ra\O_\cX^\times$ is
injective.

To see that $\O_\cX^\times\ra i_*(\O_\cD^\times)$ is surjective,
fix a point $x\in |\cD|$ and a germ
$t\in\O_{\cD,\bar{x}}^\times$. Then there is a germ
$s\in\O_{\cX,\bar{x}}$ mapping to
$t$, and this germ is invertible because
$s(\bar{x})\in\kappa(\bar{x})$ is nonzero.

It remains to see that the sequence (\ref{sequence}) is
exact in the middle at a given point
$x\in|\cX|$.
This is obvious on $\cX-\cD$, so we may assume that $x\in|\cD|$, in other words,
$\Spec(\O_{\cX,\bar{x}})$ is not irreducible.
We first check that the sequence
(\ref{sequence}) is a complex at $x$. Fix a germ
$s_{\bar{x}}\in\O_{\cX,\bar{x}}^\times$ coming from a germ
$t_{\bar{x}}\in\shB_{\cX,\bar{x}}$.
Choose an affine \'etale neighborhood $\cU\ra \cX$ so that
$s_{\bar{x}},t_{\bar{x}}$ admit representants $s,t$, and that the
canonical map
$\Spec(\O_{\cX,\bar{x}})\ra \cU$ induces a bijection on
the set of irreducible
components. Decompose $\cU=\cU_1\cup\ldots\cup \cU_n$, $n\geq 2$ into
irreducible components. Using
$$
\shB_\cU=
\shHom(\bigoplus_{i=1}^nf_*(\ZZ_{\cU_i}),\O_\cU^\times)=
\bigoplus_{i=1}^n\shHom(\ZZ_{\cU_i},\O_\cU^\times),
$$
we obtain a decomposition $t=(t_1,\ldots,t_n)$ with 
$t_i\in\Hom(\ZZ_{\cU_i},\O_\cU^\times)$, and in
turn a factorization $s=s_1\ldots s_n$ with $s_i=t_i(\rho)$.
Let $\eta_i\in \cU_i$ be the generic points.
Then $(t_i)_{\bar{\eta}_j}=1$ for $i\neq j$ because
$\shHom(\ZZ_{\cU_i},\O_\cU^\times)$ has support on $\cU_i$.
Consequently $(s_i)_{\bar{\eta}_j}=1$, and therefore
$s_i|_{\cU_j}=1$, since the $\cU_j$ have no embedded components.
Making a cyclic permutation, we calculate
$$
s_\cD = (s_1|_{\cU_2})_\cD\cdot (s_2|_{\cU_3})_\cD\cdot\ldots\cdot
(s_{n-1}|_{\cU_n})_\cD\cdot (s_n|_{\cU_1})_\cD=1.
$$
Hence $s_x$ maps to $1\in\O_{\cD,\bar{x}}^\times$, and the
sequence (\ref{sequence}) is a complex.

Finally, suppose a germ $s_{\bar{x}}\in\O_{\cX,\bar{x}}^\times$
maps to
$1\in\O_{\cD,\bar{x}}^\times$. As above, we choose an affine
\'etale neighborhood $\cU\ra \cX$ such that
$\Spec(\O_{\cX,\bar{x}})\ra \cU$ induces a bijection on the set of
irreducible components and that
$s_{\bar{x}}$ admits a representant $s$. Write $\cU=\Spec(A)$,
$\cU_i=\Spec(A_i)$, and let $t_i\in A_i$
be the image of $s\in A$. Set $B=A_1\times\ldots\times A_n$.
Then $S_\cU=\Spec(B)$,
and
$t=(t_1,\ldots,t_n)\in B$ is the image of
$s\in A$. Since $s_D=1$, we also have $t_i|_{f^{-1}(\cD)}=1$. Now the
exact sequence
$$
1\lra\O_\cX^\times\lra\O_\cS^\times\oplus\O_\cD^\times\lra
\O_{f^{-1}(\cD)}^\times\lra 1
$$
implies that each pair
$(t_i,1)\in\Gamma(\cU,\O_\cS^\times\oplus\O_\cD^\times)$ comes from a
section $s_i\in\Gamma(\cU,\O_\cX^\times)$. Sending the $i$-th
standard generator of
$f_*(\ZZ_\cS)_\cU$ to
$s_i$, we obtain a homomorphism
$h:f_*(\ZZ_\cS)_\cU\ra\O_\cU^\times$ with $h(\rho)=s$ at the generic
points. Since
$\cX$ has no embedded points, $h(\rho)=s$ holds globally. In other words
the germ
$s_x$ lies in the image of
$\rho^*:\shB_\cX\ra\O_\cX^\times$.
\qed

\medskip
To apply this calculation to log atlases we first need a comparison
result for constructible sheaves.

\begin{proposition}
\mylabel{canonical comparison}
Suppose $\cX$ is a noetherian algebraic space
satisfying Serre's condition $(S_2)$, and let $i:\cU\ra \cX$
be an open embedding containing all points
of codimension $\leq 1$. Let $\shF_1,\shF_2$ be two constructible
abelian sheaves on $\cX$. If $i_*i^{-1}(\shF_1)$ is constructible and
$i^{-1}(\shF_1)\simeq i^{-1}(\shF_2)$, then
$\shHom(\shF_1,\O_\cX^\times)\simeq \shHom(\shF_2,\O_\cX^\times)$.
\end{proposition}

\proof
Let $\shK_j,\shC_j$ be kernel and cokernel of the adjunction
maps $\shF_j\ra i_*i^{-1}(\shF_j)$, respectively. These are constructible
abelian sheaves supported by
$\cX-\cU$. Applying the functor $\shHom(.,\O_\cX^\times)$ to the exact sequences
of constructible abelian sheaves
$$
0\ra\shK_j\ra\shF_j\ra\shF_j/\shK_j\ra 0
\quadand
0\ra\shF_j/\shK_j\ra i_*i^{-1}(\shF_j)\ra \shC_j\ra 0,
$$
we reduce our problem to the following special cases:
We have a map $\shF_1\ra\shF_2$ that is either injective or surjective,
and furthermore bijective on $\cU$.

First, consider the case that we have a surjective mapping
$\shF_1\ra\shF_2$, and let $0\ra \shK\ra\shF_1\ra\shF_2\ra 0$ be the
corresponding exact sequence.
This gives an exact sequence
$$
1\lra \shHom(\shF_2,\O_\cX^\times)\lra \shHom(\shF_1,\O_\cX^\times)
\lra \shHom(\shK,\O_\cX^\times).
$$
The adjunction map $\O_\cX^\times\ra i_*i^{-1}(\O_\cX^\times)$
is injective, because $\cX$ has no embedded components,
hence there is an injection
$\shHom(\shK,\O_\cX^\times)\subset\shHom(\shK,i_*i^{-1}(\O_\cX^\times))$. We have
$$
\shHom(\shK,i_*i^{-1}(\O_\cX^\times))=
i_*\shHom(i^{-1}(\shK),i^{-1}(\O_\cX^\times))
$$
by \cite{SGA 2}, Expos\'e I, Corollary 1.5, and conclude
that $\shHom(\shF_2,\O_\cX^\times)\ra\shHom(\shF_1,\O_\cX^\times)$ is
bijective.

Second, suppose we have an injection $\shF_1\ra\shF_2$,
and let $0\ra\shF_1\ra\shF_2\ra \shC\ra 0$ be the
corresponding exact sequence. As above, we have
$\shHom(\shC,\O_\cX^\times)=1$ and obtain an exact sequence
$$
1\lra\shHom(\shF_2,\O_\cX^\times)\lra\shHom(\shF_1,\O_\cX^\times)\lra
\shExt^1(\shC,\O_\cX^\times).
$$
We finish the proof by checking that $\shExt^1(\shC,\O_\cX^\times)$ vanishes.
This is a local problem, so we may assume that $\cX$ is an affine
scheme. Let $i:\cX-\cU\ra \cX$ be the embedding of the closed subset $\cX-\cU$
of codimension $\geq 2$. Then $\shC=\shB_\cX$ for the constructible
sheaf $\shB=i^{-1}(\shC)$ on $\cX-\cU$, where $\shB_\cX=i_!(\shB)$ denotes
extension by zero.

According to \cite{SGA 4c}, Expos\'e IX, Lemma 2.10,
there are finitely many \'etale neighborhoods $\cC_i\ra \cX-\cU$, $1\leq i\leq n$
and local sections $s_i\in\Gamma(\cC_i,\shB)$ so that
the corresponding map $\bigoplus_{i=1}^n \ZZ_{\cC_i}\ra\shB$ is surjective.
We then say that $\shB$ is \emph{generated by $n$ local sections}.
Let $\shB_1\subset\shB$ be the subsheaf generated by $\ZZ_{\cC_1}$.
Using the exact sequence
$$
\shExt^1((\shB/\shB_1)_\cX,\O_\cX^\times)
\lra\shExt^1(\shB_\cX,\O_\cX^\times)\lra\shExt^1((\shB_1)_\cX,\O_\cX^\times)
$$
and induction on the number $n$ of local sections, it suffices
to treat the case that $\shB$ is generated by a single local section.
In other words, there is an exact sequence $0\ra\shB'\ra\ZZ_\cC\ra \shB\ra 0$,
where $\cC\ra \cX-\cU$ is an
\'etale neighborhood. Then we have an exact sequence
$$
\shHom(\shB'_\cX,\O_\cX^\times)\lra
\shExt^1(\shB_\cX,\O_\cX^\times)\lra\shExt^1(\ZZ_{\cC,\cX},\O_\cX^\times).
$$ 
The term on the left vanishes, and we are reduced to the case
$\shC=\ZZ_{\cC,\cX}$.

Next, choose an affine open covering $\cC_i\subset \cC$, say $1\leq i\leq m$,
so that there are affine \'etale coverings $\cV_i\ra \cX$
with $\cC_i=(\cX-\cU)\times_\cX \cV_i$
(use \cite{SGA 1}, Expos\'e I, Proposition 8.1).
Using the surjection $\bigoplus_{i=1}^m\ZZ_{\cC_i}\ra\ZZ_\cC$ and repeating the
argument in the preceding paragraph, we reduce to the case $m=1$, and write
$\cC=\cC_1$ and $\cV=\cV_1$.

Now $\cV\ra \cX$ is \'etale and $\cC\ra \cV$ is a closed embedding. Let $\cV{}'=\cV-\cC$ be
the complementary open subset.
Then we have an exact sequence of sheaves
$0\ra\ZZ_{\cV{}',\cV}\ra\ZZ_{\cV}\ra\ZZ_{\cC,\cV}\ra 0$ on $\cV$.
Extending by zero, we obtain an exact sequence
$0\ra\ZZ_{\cV{}',\cX}\ra\ZZ_{\cV,\cX}\ra\ZZ_{\cC,\cX}\ra 0$ on $\cX$,
and in turn a long exact sequence
$$
\shH^0_{\cV{}'}(\O_\cX^\times)\lra \shH^0_{\cV}(\O_\cX^\times)\lra
\shExt^1(\ZZ_{\cC,\cX},\O_\cX^\times)\lra
\shH^1_{\cV{}'}(\O_\cX^\times)\lra\shH^1_{\cV}(\O_\cX^\times).
$$
Here we applied the functor $\shExt^n(\cdot,\O_\cX^\times)$ and identified
$\shExt^n(\ZZ_{\cV,\cX},\O_\cX^\times)$ with the sheaf of local cohomology groups
$\shH^n_{\cV}(\O_\cX^\times)$ as in
\cite{SGA 2}, Expos\'e I, Proposition 2.3.
The map $\shH^0_{\cV{}'}(\O_\cX)\ra \shH^0_{\cV}(\O_\cX)$ is
surjective, because $\cX$ satisfies Serre's condition
$(S_2)$ and $\cC=\cV-\cV{}'$ has codimension $\geq 2$.
Hence, by Krull's Principal Ideal Theorem,
the map $\shH^0_{\cV{}'}(\O_\cX^\times)\ra \shH^0_{\cV}(\O_\cX^\times)$
is surjective as well.

The sheaf $\shH^1_\cV(\O_\cX^\times)$ is associated to the presheaf
$\cW\mapsto\Pic(\cV\times_\cX \cW)$.
The restriction map $\Pic(\cV\times_\cX \cW)\ra \Pic(\cV{}'\times_\cX \cW)$ is injective by
\cite{Schroeer 2001}, Lemma 1.1, so the map on sheaves
$\shH^1_{\cV{}'}(\O_\cX^\times)\ra\shH^1_{\cV}(\O_\cX^\times)$ is injective as well.
It follows that $\shExt^1(\ZZ_{\cC,\cX},\O_\cX^\times)$ vanishes as desired.
\qed

\medskip
We now apply this to our sheaf $\shA_\cX=\shHom(\sqM_\cX,\O_\cX^\times)$
of automorphism of log structures.
Let $\cS$ be the disjoint union of the irreducible components
of $\cX$ and $f:\cS\ra \cX$ the corresponding finite birational map, which is  the normalization map.

\begin{theorem}
\mylabel{canonical bijection}
Let $\cX$ be a noetherian algebraic space satisfying
Serre's condition $(S_2)$ and whose integral components of
$\Spec(\O_{\cX,\bar{x}})$ are normal
for all $x\in|\cX|$, and let
$\sqM_\cX$ be a constructible monoid sheaf with integral stalks.
Suppose there is an open subset $\cU\subset \cX$ containing all
points of codimension $\leq 1$ with $\sqM_\cU\simeq f_*(\NN_\cS)_\cU$.
Then $\shA_\cX=\shHom(f_*(\ZZ_\cS),\O_\cX^\times)$,
and we have an exact sequence
$$
1\lra \shA_\cX\lra \O_\cX^\times\lra\O_\cD^\times\lra 1,
$$
where $\cD\subset \cX$ is the branch space for the finite birational morphism
$f:\cS\ra \cX$.
\end{theorem}

\proof
To check the first assertion we apply Proposition \ref{canonical comparison}
with the constructible abelian sheaves
$\shF_1=f_*(\ZZ_\cS)$ and $\shF_2=\sqM_\cX^\gp$.
We have to check that $i_*i^{-1}f_*(\ZZ_\cS)$ is constructible,
where $i:\cU\ra \cX$ is the canonical open embedding.
We do this by showing that the adjunction map
$f_*(\ZZ_\cS)\ra i_*i^{-1}f_*(\ZZ_\cS)$ is bijective.
Fix a point $x\in |\cX|$. Then the stalks of both sides at $x$ are the free
group generated by the irreducible components of $\Spec(\O_{\cX,\bar{x}})$,
and bijectivity follows.

Having $\shA_\cX=\shHom(f_*(\ZZ_\cS),\O_\cX^\times)$, the second assertion
directly follows
from Proposition \ref{exactness}.
\qed

\section{The restricted conormal sheaf}
\mylabel{restricted conormal sheaf}

We now use the exact sequence from Theorem
\ref{canonical bijection} to compute the cohomology group $H^2(\cX,\shA_\cX)$,
which contains the obstruction for the existence of a global log structure.
We also compute the cohomology group $H^1(\cX,\shA_\cX)$, which measures how many
isomorphism classes of global log structures exists.
Throughout, we make the following assumptions:
Let $\cX$ be a reduced noetherian algebraic space satisfying Serre's
condition $(S_2)$ and such that for all $x\in|\cX|$ the integral components
of $\Spec(\O_{\cX,\bar{x}})$ are normal. Furthermore,  $\sqM_\cX$ is a
constructible monoid sheaf with  integral stalks satisfying the
conditions of Theorem  \ref{canonical bijection}.
We set $\shA_\cX=\shHom(\sqM_\cX,\O_\cX^\times)$.

Consider the short exact sequence
$$
1\lra \shA_\cX\lra\O_\cX^\times\lra
i_*(\O_\cD^\times)\lra 1,
$$
where $i:\cD\ra \cX$ is the closed embedding of the space of nonnormality.
We have $\Pic(\cD)=H^1(\cX,i_*\O_\cD^\times)$, because
$R^1i_*(\O_\cD^\times)=0$ by Hilbert's Theorem 90. The preceding short exact
sequence gives a long exact sequence
\begin{equation}\label{long exact sequence}
\Pic(\cX)\lra\Pic(\cD)\lra H^2(\cX,\shA_\cX) \lra \Br'(\cX),
\end{equation}
where $\Br'(\cX)=H^2(\cX,\O_\cX^\times)$ is the cohomological Brauer group.
We see that an $\shA_\cX$-gerbe $\shG$ faces two obstructions against
neutrality: The first obstruction is the image of the gerbe class
$[\shG]\in H^2(\cX,\shA_\cX) $ in the cohomological Brauer group
$\Br'(\cX)$. This obstruction vanishes if and only if there is an
invertible $\O_\cD$-module  $\shN_\cD$ whose $\shA_\cX$-gerbe of extensions
to invertible $\O_\cX$-modules is equivalent to $\shG$. Once the first
obstruction vanishes, the second obstruction is the extendibility of
$\shN_\cD$ to $\cX$.  It turns out that, under suitable assumptions, the Brauer obstruction
vanishes automatically:

\begin{proposition}
\mylabel{brauer obstruction}
Suppose there is a global section $\rho\in\Gamma(\cX,\sqM_\cX)$ such that
the stalks $\rho_{\bar\eta}$ generate $\sqM_{\cX,\bar{\eta}}=\NN$ for all generic
points
$\eta\in|\cX|$. Let $\shG$ be a log atlas on $\cX$ with respect to $\sqM_\cX$.
Then the gerbe class $[\shG]\in H^2(\cX,\shA_\cX)$ maps to zero in
the cohomological Brauer group $\Br'(\cX)$.
\end{proposition}

\proof
First note that the map $\shA_\cX\ra\O_\cX^\times$ from
Theorem \ref{canonical bijection} is nothing but the
evaluation map $\rho^*(h)=h(\rho)$.

Let $\shP$ be the gerbe of invertible sheaves on
\'etale neighborhoods $\cU\ra \cX$. This $\O_\cX^\times$-gerbe represents
the zero element in $H^2(\cX,\O_\cX^\times)$.
To check that $[\shG]$ maps to zero in $\Br'(\cX)$, we have to
construct a cartesian functor
$\shG\ra\shP$ that is equivariant with respect to the map 
$\rho^*:\shA_\cX\ra\O_\cX^\times$,
as explained in \cite{Giraud 1971}, Chapter IV, Definition 3.1.4.
Let $U\in\shG$ be a log space. The exact sequence
$$
1\lra \O_U^\times\lra \shM_U^\gp\lra \sqM_U^\gp\lra 0
$$
yields a coboundary map $H^0(\cU,\sqM_U^\gp)\ra\Pic(\cU)$.
Let $\shN_U$ be the invertible $\O_\cU$-module associated
to the $\O_\cU^\times$-torsor
$\shM_U\times_{\sqM_\cU}\left\{\rho_U\right\}$.
Then $U\mapsto \shN_U$
is the desired cartesian functor $\shG\ra\shP$. You easily check that the diagram
$$
\begin{CD}
H^0(\cU,\shA_\cX) @>\rho^*>> H^0(\cU,\O_\cX^\times)\\
@V\simeq VV@VV\simeq V\\
\Aut(U/\cU) @>>>\Aut(\shN_U)
\end{CD}
$$
is commutative and compatible with restrictions.
This means  $\rho^*([\shG])=[\shP]=0$.
\qed

\medskip
From now on we assume that a section $\rho\in\Gamma(\cX,\sqM_\cX)$ as in
Proposition
\ref{brauer obstruction} exists.
Then we see that the gerbe class $[\shG]\in H^2(\cX,\shA_\cX)$ of a log atlas
$\shG$ comes from an invertible $\O_\cD$-module. It turns out
that there is a \emph{canonical choice}
as follows:
Pick an \'etale covering $\cU\ra \cX$ admitting a log space
$U\in\shG$. Passing to a finer covering, we also have a section
$\tilde{\rho}\in\Gamma(\cU,\shM_U)$ mapping to
$\rho\in\Gamma(\cU,\sqM_\cX)$. Next choose an \'etale covering
$\cV\ra \cU\times_\cX \cU$ so that there is an isomorphism
$\phi:p_1^*(U)\ra p_0^*(U)$, which is given by a
bijection
$\phi:p_1^*(\shM_U)\ra p_0^*(\shM_U)$
fixing the subsheaf
$\O_\cV^\times$ pointwise and inducing the identity on the quotient sheaf
$\sqM_\cV$. As explained in Section \ref{cech cohomology}, the equation
$$
c\cdot p_1^*(\phi) = p_0^*(\phi)p_2^*(\phi)\in
\Isom(p_1^*p_1^*U, p_0^*p_0^*U)
$$
defines a 2-cocycle $c\in Z^2(\cV/\cU,\shA_\cX)$ representing the gerbe
class of the log atlas
$\shG$.
Now the equation
\begin{equation}
\label{cochain}
\phi(p_0^*(\tilde{\rho})) = e\cdot p_1^*(\tilde{\rho})\in
\Gamma(\cV,p_1^*(\shM_U))
\end{equation}
defines a cochain $e\in C^1(\cV/\cU,\O_\cX^\times)$. We claim that its
restriction
$e_\cD$ to the ramification locus $\cD\subset \cX$ of $f:\cS\ra \cX$ becomes a
cocycle. Indeed, using the simplicial identities
$p_j^*p_i^*=p_i^*p_{j-1}^*$, $i<j$, we compute
\begin{equation}
\label{comparison}
\begin{split}
p_1^*(\phi)(p_0^*p_0^*(\tilde{\rho})) &=
p_1^*p_1^*(\tilde{\rho})\cdot p_1^*(e),\\
p_0^*(\phi)p_2^*(\phi)(p_0^*p_0^*(\tilde{\rho}))&=
p_1^*p_1^*(\tilde{\rho})\cdot p_0^*(e)p_2^*(e).
\end{split}
\end{equation}
On the other hand, the two isomorphisms $p_1^*(\phi)$ and
$p_2^*(\phi) p_0^*(\phi)$ differ by $c$, and $c(\rho)_\cD=1$
according to Proposition \ref{exactness}, hence
$p_1^*(e_\cD)=p_0^*(e_\cD)p_2^*(e_\cD)$.

The cocycle $e_\cD\in Z^1(\cV/\cU,\O_\cD^\times)$ defines an invertible
$\O_\cD$-module
$\shN_\cD$. In fact, its isomorphism class is an invariant
of the log atlas $\shG$:

\begin{proposition}
\mylabel{conormal sheaf}
The isomorphism class of 
$\shN_\cD$ depends only
on the log atlas $\shG$ and the section $\rho\in\Gamma(\cX,\sqM_\cX)$. It maps to
the gerbe class
$[\shG]$ under the coboundary map $\Pic(\cD)\ra H^2(\cX,\shA_\cX)$.
\end{proposition}

\proof
Replacing the \'etale coverings $\cV\rightrightarrows \cU$ by some
refinement replaces the cocycle $e_\cD\in Z^1(\cV/\cU,\O_\cD^\times)$ by its
restriction to some finer covering.
Changing the lift $\tilde{\rho}\in \Gamma(U,\shM_U)$ by some
invertible function changes the cocycle
$e_\cD$ by a coboundary.
Modifying the bijection $\phi:p_1^*(\shM_U)\ra p_0^*(\shM_U)$ with 
some 
$h\in C^1(\cV/\cU,\shA_\cX)$ does not affect $e_\cD$ at all,
because $h(\rho)_\cD=1$
by Proposition
\ref{exactness}. Summing up, the isomorphism class
$\shN_\cD\in\Pic(\cD)$ does not depend on our choices.

To calculate the coboundary $\partial(\shN_\cD)$, we use the cochain
$e\in C^1(\cV/\cU,\O_\cX^\times)$
from Equation (\ref{cochain}) as a lift for the cocycle $e_\cD$.
Then the cocycle $h\in Z^2(\cV/\cU,\shA_\cX)$ defined by
$h(\rho)=p_0^*(e)p_2^*(e)/p_1^*(e)$ represents
$\partial(\shN_\cD)$. On the other hand, the bijections
$p_1^*(\phi)$ and $p_0^*(\phi)p_2^*(\phi)$ differ
by $p_0^*(e)p_2^*(e)/p_1^*(e)$ on
$p_1^*p_1^*(\tilde{\rho})=p_2^*p_1^*(\tilde{\rho})$, according to
Equation (\ref{comparison}). By Proposition \ref{exactness},
this means that these bijections differ by $h$, and we conclude
$\partial(\shN_\cD)=[\shG]$.
\qed

\medskip
By abuse of notation, we call the invertible $\O_\cD$-module
$\shN_\cD$ in Proposition \ref{conormal sheaf} the
\emph{restricted conormal sheaf} of the log atlas $\shG$.
The main result of this section is the following classification
result:

\begin{theorem}
\mylabel{existence}
There is a global log structure 
$X\in\shG$ if and only if the restricted
conormal sheaf $\shN_\cD$ extends to an invertible $\O_\cX$-module.
\end{theorem}

\proof
Proposition \ref{exactness} gives an exact sequence 
$$
\Pic(\cX)\lra\Pic(\cD)\lra H^2(\cX,\shA_\cX).
$$
According to Proposition \ref{conormal sheaf}, the restricted
conormal sheaf
$\shN_\cD$ maps to the gerbe class of $\shG$, and the assertion
follows.
\qed 

\medskip
For the rest of this section we study the action of
$H^1(\cX,\shA_\cX)$ on the isomorphism class of global log structures
$X\in\shG$. First note that
each $U\in\shG$ comes along with an exact sequence of
abelian groups
$$
1\lra\O_\cU^\times\lra\shM^\gp_U\lra\sqM^\gp_U\lra 0,
$$
and defines a $\O_\cU^\times$-torsor $\shM^\gp_U\times_{\sqM^\gp_U}
\left\{\rho_U\right\}$, hence an invertible
$\O_\cU$-module
$\shN_U$. We call $\shN_U$ the
\emph{conormal sheaf} of the log structure. Its restriction to $\cD$
is isomorphic to $\shN_\cD$, by the very definition of the restricted
conormal sheaf below Equation
(\ref{cochain}).

\begin{proposition}
\mylabel{induced}
Let $X\in\shG$ be a global log structure, $\shN_X$ its conormal
sheaf,
$\alpha\in H^1(\cX,\shA_\cX)$ a cohomology class, and
$\shL=\rho^*(\alpha)$ its
image in
$\Pic(\cX)$. Then the conormal sheaf of the global log
structure $X+\alpha\in\shG$ is isomorphic to
$\shN_X\otimes\shL$.
\end{proposition}

\proof
Choose an \'etale covering $\cU\ra \cX$ and a cocycle $h\in
Z^1(\cV/\cU,\shA_\cX)$ representing the cohomology class
$\alpha$. Here $\cV=\cU\times_\cX \cU$.
Let $\phi:p_1^*(U)\ra p_0^*(U)$ be the canonical
isomorphism such that $(U,\phi)$ is a descent datum for
$X$. Consequently $(U,\phi h)$ is a descent datum for
$X+\alpha$. Refining $\cU$, we may also choose a lift
$\tilde{\rho}\in\Gamma(\cU,\shM_U)$ for
$\rho$. Then the cocycle $e\in Z^1(\cV/\cU,\O_\cX^\times)$ defined by
$e\cdot p_1^*(\tilde{\rho})=p_0^*(\tilde{\rho})$ represents the
conormal sheaf
$\shN_X$. It follows that $e\cdot h(\rho)$ is both a cocycle for the
conormal sheaf of $X+\alpha$ and the tensor product
$\shN_X\otimes\shL$.
\qed

\begin{corollary}
Let $\shN$ be an invertible $\O_\cX$-module extending
the restricted conormal sheaf $\shN_\cD$.
Then the set of isomorphism classes of global log
spaces $X\in\shG$ whose conormal sheaf is isomorphic to
$\shN$
is a torsor for the cokernel $\Gamma(\O_\cD^\times)/\Gamma(\O_\cX^\times)$
of the restriction map $\Gamma(\O_\cX^\times)\ra\Gamma(\O_\cD^\times)$.
\end{corollary}

\proof
There is a global log space $X\in\shG$ by Theorem
\ref{existence}, and its conormal sheaf $\shN_X$ extends the
restricted conormal sheaf
$\shN_\cD$. Proposition \ref{exactness} gives an exact sequence
$$
\Gamma(\O_\cX^\times)\lra \Gamma(\O_\cD^\times) \lra
H^1(\cX,\shA_\cX)\lra \Pic(\cX)\lra\Pic(\cD).
$$
Using Proposition \ref{induced}, we
we may change the global log structure
$X$ by some element in $H^1(\cX,\shA_\cX)$ so that
its conormal sheaf becomes isomorphic to
$\shN$. Moreover, all such log structures differ by elements in the
subgroup
$\Gamma(\O_\cD^\times)/\Gamma(\O_\cX^\times)$.
\qed

\medskip
We can say more about the action of the
subgroup $\Gamma(\O_\cD^\times)/\Gamma(\O_\cX^\times)\subset
H^1(\cX,\shA_\cX)$ on the global log structures:

\begin{proposition}
The action of $\Gamma(\O_\cD^\times)/\Gamma(\O_\cX^\times)$ on the set
of isomorphism classes of global log structures $X\in\shG$
does not change the sheaf of sets $\shM_X$ and the surjective
map $\shM_X\ra\sqM_\cX$.
\end{proposition}

\proof
Given an invertible function $s_\cD\in H^0(\cD,\O_\cD^\times)$, choose an
\'etale covering $\cU\ra \cX$ so that $s_\cD$ extends to a cochain $s\in
\Gamma(\cU,\O_\cX^\times)$. Then there is a 1-cocycle $h\in
Z^1(\cV/\cU,\shA_\cX)$ with $h(\rho)=p_0^*(s)/p_1^*(s)$, where
$\cV=\cU\times_\cX\cU$.

Given a log space $X\in\shG$, the canonical isomorphism
$\phi:p_0^*(U)\ra p_1^*(U)$ yields a descent datum
$(U,\phi)$ defining the log space $X$. As discussed
before Remark \ref{Cech or sheaf},  $(U,\phi h)$ is another
descent datum defining another log space
$X'=(\cX,\shM_{X'},\alpha_{X'})$, and the torsor $\Isom(X,X')$
corresponds to the cohomology class of the coboundary
$\partial(s_\cD)\in H^1(\cX,\shA_\cX)$.  
We now exploit that the 1-cocycle $h$ is defined in terms
of $s$, which lives inside $\O_\cX^\times\subset\shM_X$:
Indeed, the commutative diagram
$$
\begin{CD}
p_1^*(\shM_U) @>\phi h>> p_0^*(\shM_U)\\
@Vp_1^*(s) VV @VVp_0^*(s)V\\
p_1^*(\shM_U) @>>\phi >p_0^*(\shM_U)
\end{CD}
$$
constitutes a bijection of descent data, hence
a bijection of set-valued sheaves $\shM_{X'}\ra\shM_X$.
This map is compatible with the surjections to $\sqM_\cX$,
because the images of $p_i^*(s)$ in $\sqM_\cV$ vanish.
\qed

\section{Gorenstein toric varieties}

Our next goal is to study log atlases whose log spaces
$U\in\shG$ are locally isomorphic to a boundary divisors in
toroidal embeddings.  We come to this in the next section. Here we
collect some facts on boundary divisors in toric varieties, which we
shall use later.

Fix a ground field $k$ of arbitrary characteristic $p\geq 0$. Recall
that affine \emph{toric varieties} are of the form 
$\cZ=\Spec k[\sigma^\vee\cap M]$. Here $M$ is a finitely generated free abelian
group, $\sigma$ is a convex rational polyhedral cone in
$N\otimes_\ZZ\RR$ not containing nontrivial linear subspaces, and
$N=\Hom(M,\ZZ)$. Note that monoids of the form $P=\sigma^\vee\cap M$
are precisely the fine saturated torsionfree monoids, and we have
$M=P^\gp$. Here \emph{saturated} means that that each $p\in P^\gp$
with $np\in P$ for some integer $n>0$ lies in $P$.

From now on we usually write $P=\sigma^\vee\cap M$. To avoid
confusion  of the additive composition law for the monoid $P$ and 
the multiplicative composition law for the ring $k[P]$, we use
exponential notation $\chi^p\in k[P]$ for elements $p\in P$. We refer
to the books of Kempf et al.\ \cite{Kempf et al. 1973} and Oda
\cite{Oda 1988} for the theory or toric varieties and toroidal
embeddings.

The inclusion of monoids $P\subset k[P]$ defines a log space
$Z$ with underlying space $\cZ$. Its ghost sheaf $\sqM_Z=\shM_Z/\O_Z^\times$ is nothing but
the sheaf of effective Cartier divisors that are invariant under the canonical action of
the torus $\cT=\Spec k[M]$.  Consider the complement $\cZ_0=\cZ-\cT$ endowed
with its reduced structure. We call $\cZ_0$ the \emph{boundary divisor}
of the affine toric variety $\cZ$. It inherits the structure of a log space
$Z_0$ from the ambient log space $Z$.
From now on we denote by $Z,Z_0$ the log spaces whose underlying schemes
are  toric varieties and their boundary divisors, respectively.

The reflexive rank one sheaf $\O_\cZ(\cZ_0)$ corresponding to the
boundary divisor $\cZ_0\subset \cZ$ is a dualizing sheaf for $\cZ$,
according to \cite{Oda 1988}, Corollary 3.3 and the Remark thereafter.
Consequently, the Weil divisor  $\cZ_0\subset \cZ$ is Cartier if and only if the
toric variety $\cZ$ is Gorenstein.  In terms of the cone $\sigma$, this
means that there is an element $\rho_\sigma\in\sigma^\vee\cap M$ such
that the linear form $\rho_\sigma\in N^\vee$ takes value $1$ on the
integral generator of each 1-dimensional face
$\sigma_i\subset\sigma$. In terms of the monoid $P=\sigma^\vee\cap M$, this
translates into  the following condition:
There is a unique element $\rho_\sigma\in P$ with $\rho_\sigma+P = \interiour(P)$,
as Stanley explained in \cite{Stanley 1978}, Theorem 6.7.
Here $\interiour(P)=(\interiour\sigma^\vee)\cap M$ is the set of lattice points inside
the topological interiour $\interiour\sigma^\vee$ of the real cone $\sigma^\vee$.

We are mainly interested in this situation. Then the Cartier divisor
$\cZ_0\subset \cZ$ corresponds to the section $\rho_\sigma\in\Gamma(\cZ,\sqM_\cZ)=P$,
and we shall also denote by $\rho_\sigma\in\Gamma(\cZ_0,\sqM_{\cZ_0})=P$ the
induced section. To summarize the situation:

\begin{proposition}
\mylabel{local sing}
Let $\cZ=\Spec(k[\sigma^\vee\cap M])$ be a Gorenstein toric variety.
Then the boundary divisor $\cZ_0$ is Cohen-Macaulay,
Gorenstein, generically smooth, and has normal crossing singularities in
codimension one. 
\end{proposition}

\proof
Without any hypothesis, the schemes $\cZ$ and $\cZ_0$ are
Cohen-Macaulay by Ishida's Criterion (see \cite{Oda 1988}, page 126).
By assumption, $\cZ$ is Gorenstein and $\cZ_0$ is Cartier, so $\cZ_0$
is Gorenstein as well.
The toric variety $\cZ$ is smooth in codimension $\leq 1$, and 
has $A_n$-singularities in codimension two.
Saying that a point $z\in \cZ$ of codimension two has
an $A_n$-singularity means that the complete local ring
$\O_{\cZ,z}^\wedge$ is isomorphic to $\kappa(z)[[x^{n+1},y^{n+1},xy]]$.
A local computation shows that $\cZ_0$ is generically smooth,
and is a Cartier divisor inside an $A_n$-singularity in codimension one.
Therefore $\cZ_0$ has normal crossing singularities
in codimension one.
\qed

\medskip
Let us now consider the ghost sheaf $\sqM_{Z_0}$ of the log space $Z_0$. Later, 
we have to glue isomorphic copies of such sheaves.
The following result tells us that the cocycle condition
then holds automatically:

\begin{proposition}
\mylabel{local ghost}
Let $Z=\Spec(k[\sigma^\vee\cap M])$ be a toric variety with its canonical log structure. Then the sheaf of
groups $\shAut(\sqM_{Z_0})$ is trivial.
\end{proposition}

\proof
We have to check that $\shAut(\sqM_{Z_0})_{\bar{x}}=0$ for a given point
$x\in Z_0$. Clearly we may assume that $x$ lies in the closed orbit.
The generic points $\eta_i\in Z_0$ correspond to the invariant Weil divisors
on $\cZ$, which correspond to the extremal rays $\sigma_i\subset\sigma$.
We have $\sqM_{Z_0,\bar{x}}=\sigma^\vee\cap M/\sigma^\perp\cap M$, and
the localization map $\sqM_{Z_0,\bar{x}}\ra\sqM_{Z_0,\bar{\eta}_i}$ is
nothing but the canonical map to $\sigma_i^\vee\cap M/\sigma_i^\perp\cap
M$. The direct sum of these maps
$$
\sigma^\vee\cap M/\sigma^\perp\cap M \lra\bigoplus_i
\sigma_i^\vee\cap M/\sigma_i^\perp\cap M
$$
is injective. Since
any automorphism of $\sqM_{Z_0}$ obviously induces the identity on
$\sqM_{Z_0,\bar{\eta_i}}=\NN$, it has to induce the identity on
$\sqM_{Z_0,\bar{x}}$ as well.
\qed

\medskip
We now turn to 
a problem that occurs if $\cZ$ is singular in codimension two:
Although $\cZ_0$ is normal crossing in codimension one, the ghost sheaf
$\sqM_{Z_0}$ does not look like the ghost sheaf of a
normal crossing singularity. 
But we definitely need this property to
apply Theorem \ref{canonical bijection}.
To overcome this problem we 
make another assumption, namely that the toric variety $\cZ$ satisfies the
regularity condition
$(R_2)$, in other words, $\cZ$ is regular in codimension $\leq 2$.
In terms of the cone $\sigma\subset N\otimes\RR$, this means that
for each 2-dimensional face $\sigma'\subset\sigma$, the two integral vectors
generating $\sigma'$ form a basis for
$(\sigma'-\sigma')\cap N$, which is a free abelian group of rank two.

Let $\cS$ be the disjoint union of the irreducible components of $\cZ_0$,
and $f:\cS\ra \cZ_0$ the canonical map. Note that this is in fact the
normalization of
$\cZ_0$.

\begin{proposition}
\mylabel{local irred}
Let $Z=\Spec(k[\sigma^\vee\cap M])$ be a toric variety
satisfying regularity condition $(R_2)$, endowed with its canonical log structure.
Then there is an open subset $U\subset Z_0$ containing all points
of codimension $\leq 1$ such that $\sqM_U\simeq f_*(\NN_S)|_U$.
\end{proposition}

\proof
This is a local problem because $\shAut(\sqM_{Z_0})=0$ by Proposition
\ref{local ghost}.
Replacing $Z$ by some affine invariant open subsets, we may
assume that $Z$ is regular, and then the assertion is trivial.
\qed

\medskip
Summing up, we can say that for boundary divisors $\cZ_0$
in Gorenstein toric varieties $\cZ=\Spec(k[\sigma^\vee\cap M])$ satisfying
regularity condition
$(R_2)$, our results from Section \ref{sheaf of automorphisms}
and Section \ref{restricted conormal sheaf} do apply.

\section{Gorenstein toroidal crossings}

In this section we explore log atlases whose 
log spaces $U\in\shG$ are locally boundary divisors in Gorenstein
toroidal embeddings that are regular in codimension $\leq 2$. Throughout, we
fix a ground field
$k$ of characteristic $p\geq 0$, and let $\cX$ be an algebraic $k$-space
of finite type. We also fix a constructible monoid sheaf
$\sqM_\cX$ with surjective specialization maps and 
a global section $\rho\in\Gamma(\cX,\sqM_\cX)$.

Suppose we have a log atlas $\shG$ on $\cX$ with respect to $\sqM_\cX$.
A \emph{gtc-chart} consists of the following:
A Gorenstein toric variety $Z=\Spec k[\sigma^\vee\cap M]$ viewed as a log space and satisfying the
regularity condition
$(R_2)$, an affine scheme $\cU$ endowed with \'etale maps $\cX\la \cU\ra \cZ_0$,
and a bijection $\varphi:\sqM_{Z_0}|_\cU\ra\sqM_{\cX}|_\cU$, such
that the following two conditions hold:
First, the bijection $\varphi$ maps the section
$\rho_\sigma|_\cU\in\Gamma(\cU,\sqM_{Z_0})$ corresponding to the Cartier divisor
$\cZ_0\subset \cZ$ to our given section $\rho|_\cU$.
Second, we have
$(U,\varphi)\in\shG$, where $U$ is the log structure induced from
the log space $Z_0$.

By abuse of notation, we usually omit the toric variety $Z$
and the identification
$\varphi$ from the notation and speak about gtc-charts $\cX\la \cU\ra \cZ_0$.
Moreover, we say that a given point $x\in|\cX|$ lies in a gtc-chart 
$\cX\la \cU\ra \cZ_0$ if it is in the image of $|\cU|\ra |\cX|$.

\begin{definition}
\mylabel{gtc-atlas}
A log atlas $\shG$ on $\cX$ with respect to $\sqM_\cX$ is called a
\emph{gtc-atlas} if each point $x\in|\cX|$ lies in at least one
gtc-chart $\cX\la \cU\ra \cZ_0$.
\end{definition}

The symbol \emph{gtc} abbreviates \emph{Gorenstein toroidal crossings}.
This terminology is justified as follows:
According to \cite{SGA 1}, Expos\'e I, Proposition 8.1, there
is an \'etale covering $\cU{}'\ra \cU$ and an \'etale map $\cZ{}'\ra \cZ$ fitting
into a
cartesian diagram
$$
\begin{CD}
\cU{}' @>>> \cZ{}'\\
@VVV @VVV\\
\cZ_0 @>>> \cZ.
\end{CD}
$$
Note that ${\cZ}{}'-\cU{}'\subset \cZ{}'$ is a \emph{toroidal embedding} (see \cite{Kempf
et al. 1973}, Definition 1 on page 54), so gtc-charts locally identify
$\cX$ with the boundary divisor of a Gorenstein toroidal embedding.
If $\cZ$ is a regular toric variety, then $\cX$ has normal crossing
singularities. The notion of gtc-charts generalize normal crossing
singularities to a broader class of singularities, which one might
call \emph{Gorenstein toroidal crossing} singularities.

The existence of a gtc-atlas poses certain local conditions
on the algebraic space $\cX$. Let $\cS$ be the disjoint union
of the irreducible components of $\cX$ and $f:\cS\ra \cX$ the corresponding
birational finite map.

\begin{proposition}
Suppose $\cX$ admits a gtc-atlas $\shG$ with respect to $\sqM_\cX$.
Then $\cX$ is Cohen-Macaulay, Gorenstein, reduced, and
has normal crossing singularities in codimension $\leq 1$.
There is an open subset $\cU\subset \cX$ containing
all points of codimension $\leq 1$ such that $\sqM_\cU=f_*(\NN_\cS)|_\cU$.
Furthermore, $\rho_{\bar\eta}$ generates 
$\shM_{\cX,\bar\eta}=\NN$ for each generic point $\eta\in \cX$.
\end{proposition}

\proof
The first assertion follows from the
corresponding properties for boundary divisors in 
Gorenstein toric varieties satisfying regularity condition
$(R_2)$, as in Proposition \ref{local sing}.
The second assertion is local by Proposition \ref{local ghost},
and therefore follows from Proposition \ref{local irred}.
The last assertion is obvious.
\qed

\medskip
This tells us that the results from Section
\ref{sheaf of automorphisms} and Section \ref{restricted conormal sheaf}
do apply.
In particular, a gtc-atlas $\shG$ comes along with
its restricted conormal sheaf $\shN_\cD$ 
on the subspace of nonnormality $\cD\subset \cS$, and $\shG$ admits a global
log space
$X\in\shG$ if and only if the restricted conormal sheaf
extends to an invertible sheaf on $\cX$.

Our next goal is to relate gtc-atlases to local
infinitesimal deformations.
Suppose $\shG$ is a gtc-atlas on $\cX$ with respect to $\sqM_\cX$ and $\rho$.
Fix a point $x\in|\cX|$ and choose a gtc-chart $\cX\la \cU\ra \cZ_0$ containing $x$,
with $\cU$ affine.
Let $Z=\Spec k[\sigma^\vee\cap M]$ be the corresponding Gorenstein toric
variety viewed as a log space, $\rho_\sigma\in\sigma^\vee\cap M$ the monomial defining the
Cartier divisor $\cZ_0\subset \cZ$, and 
$\chi^{\rho_\sigma}\in k[\sigma^\vee\cap M]$ the corresponding equation.
Then $\chi^{2\rho_\sigma}$ defines another Cartier divisor $\cZ_1\subset \cZ$,
and $\cZ_0\subset \cZ_1$ is an infinitesimal extension with ideal
$(\chi^{\rho_\sigma})/(\chi^{2\rho_\sigma})\simeq
\O_{\cZ_0}$. According to
\cite{SGA 1}, Expos\'e I, Theorem 8.3, there is an \'etale map $\cU_1\ra \cZ_1$
fitting into a cartesian diagram
$$
\begin{CD}
\cU @>>> \cU_1\\
@VVV @VVV\\
\cZ_0 @>>> \cZ_1,
\end{CD}
$$
and $\cU\subset \cU_1$ is a first order extension with ideal
$\O_\cU$. The isomorphism class of such extensions
correspond to classes in
$$
\Ext^1(\Omega^1_{\cU/k},\O_\cU) = H^0(\cU,\shExt^1(\Omega^1_{\cX/k},\O_\cX)).
$$
The latter groups are isomorphic because we assumed that $\cU$ is affine.
Of course, the class of $\cU\subset \cU_1$ depends on the choice of the
gtc-chart
$\cX\la \cU\ra \cZ_0$ and
the \'etale map $\cU_1\ra \cZ_1$. However, we get rid of this dependence if
we pass to the limit and allow rescaling:

\begin{proposition}
\mylabel{extension class}
The $\O_{\cX,\bar{x}}$-submodule in
$\shExt^1(\Omega^1_{\cX/k},\O_\cX)_{\bar{x}}$ generated by the extension
class of
$\cU\subset \cU_1$ depends only on the gtc-atlas $\shG$.
\end{proposition}

\proof
Suppose we have two gtc-charts $\cX\la \cU\ra \cZ_0$ and $\cX\la \cU{}'\ra \cZ{}'_0$
containing
$x$, with certain affine Gorenstein toric varieties $Z=\Spec k[P]$ and 
$Z'=\Spec k[P']$.
Replacing $\cU$ and $\cU{}'$ by some common affine \'etale neighborhood, we may
assume $\cU=\cU{}'$. Choose a point $u\in \cU$ representing $x$,
let $f:\cU\ra \cZ_0$ and $f':\cU\ra \cZ{}_0'$ be the canonical maps,
and set $z=f(u)$ and $z'=f'(u)$.

Recall that among the toric orbits in the toric variety $Z$
there is a minimal toric orbit, which is the unique closed toric orbit.
Replacing $P$ by a suitable localization $P + f\ZZ$, $f\in P$,
and $\cU$ by an open subset, we may assume that the points $f(u)\in Z$ and
$f'(u)\in Z'$ are contained in  the closed toric orbit.
We then have $P/P^\times=\sqM_{X,\bar{x}}=P'/{P'}^\times$.
This identification of monoids extends to an identification of groups $(P/P^\times)^\gp=(P'/P'^\times)^\gp$,
because the monoids in question are saturated.
Moreover,
the free abelian groups $P^\times$ and ${P'}^\times$ have the same rank,
because both
$\dim(\cZ_0)$ and $\dim(\cZ{}'_0)$ equal the dimension of $\cX$ in a neighborhood
of $x$.
We infer that there is an (uncanonical) bijection
$b:P\ra P'$ covering the canonical identification
$P/P^\times=P'/{P'}^\times $. The morphism $f:\cU\ra \cZ$ is defined via the
composition
\begin{equation}
\label{defining composition}
P_\cU\lra\shM_{Z_0}|_U\lra f^*(\shM_{Z_0}) \stackrel{\alpha}{\lra}\O_\cU,
\end{equation}
and the analogous statement holds for $f':\cU\ra \cZ{}'$.
The commutative diagram
$$
\begin{CD}
\Spec(\O_{\cX,\bar{x}})
@>f'_{\bar{x}}>>\Spec(\O_{\cZ{}'_0,\bar{z}'})
@>>>\Spec(\O_{\cZ{}',\bar{z}'}) \\ @V\id VV@VVg_0 V\\
\Spec(\O_{\cX,\bar{x}})
@>>f_{\bar{x}}>\Spec(\O_{\cZ_0,\bar{z}})
@>>>\Spec(\O_{\cZ_0,\bar{z}}) 
\end{CD}
$$
defines a bijection $g_0=f_{\bar{x}}{f'}^{-1}_{\bar{x}}$. Note that
$f_{\bar{x}}$ and $f'_{\bar{x}}$ are isomorphisms, because
$f$ and $f'$ are
\'etale. We now seek to construct a bijection
$$
g:\Spec(\O_{\cZ{}'_0,\bar{z}'})\lra\Spec(\O_{\cZ_0,\bar{z}})
$$ 
extending 
$g_0$. Replacing
$\cU$ by some smaller affine \'etale neighborhood, we may assume that there is
an isomorphism of log structures $\phi:f^*(\shM_{Z_0})\ra
{f'}^*(\shM_{Z'_0})$. We have inclusions of sheaves $P_\cU\subset
f^*(\shM_{Z_0})$ and $P'_\cU\subset {f'}^*(\shM_{Z_0})$, and these constant
submonoid sheaves both surject onto $\sqM_U$. Consequently, the equation
$$
\phi(p)=h(p)\cdot b(p),\quad p\in P
$$
inside the stalk ${f'}^*(\shM_{Z_0})_{\bar{x}}$ defines a map
$h:P\ra\O_{\cX,\bar{x}}^\times$. As in the proof of Proposition
\ref{automorphism sheaf}, we infer that $h$ is a homomorphism of monoids.
Since $P^\gp$ is free, we may lift $h$ to a monoid homomorphism
$h:P\ra\O_{\cZ{}',\bar{z}'}^\times$.

To proceed, let $k[P']^\sh=\O_{Z',\bar{z}'}$ be the strict henselization
of $k[P']$ at the prime ideal corresponding to $z'\in Z'$. The map
$P\ra k[P']^\sh$, $p\mapsto h(p)\chi^{b(p)}$ defines a homomorphism
$k[P]\ra k[P']^\sh$, which by (\ref{defining composition}) makes
the diagram 
$$
\begin{CD}
\Spec(\O_{\cX,\bar{x}}) 
@>f'>>\Spec(\O_{\cZ{}',\bar{z}'}) \\
@V\id VV@VV V\\
\Spec(\O_{\cX,\bar{x}})
@>>f>\Spec(k[P]) 
\end{CD}
$$
commutative.
Therefore the preimage of the maximal ideal in $k[P']^\sh$
under the map $k[P]\ra k[P']^\sh$
is the prime ideal in $k[P]$
corresponding to 
$z\in \cZ$. In turn, we obtain a homomorphism
$k[P]^\sh\ra k[P']^\sh$, where $k[P]^\sh=\O_{\cZ,\bar{z}}$
is the strict henselization of $k[P]$ at the prime ideal for $z\in Z$.
This homomorphism defines the desired morphism $g$ making the
diagram
$$
\begin{CD}
\Spec(\O_{\cX,\bar{x}}) @>>>\Spec(\O_{\cZ_0,\bar{z}})
@>>>\Spec(\O_{\cZ,\bar{z}}) \\ 
@V\id VV @VVg_0 V @VVg V\\
\Spec(\O_{\cX,\bar{x}}) @>>>\Spec(\O_{\cZ{}'_0,\bar{z}'})
@>>>\Spec(\O_{\cZ{}',\bar{z}'})
\end{CD}
$$
commutative.
The rest is easy: Choose affine \'etale neighborhoods $\cW\ra \cZ$ and
$\cW{}'\ra \cZ{}'$ so that there is an isomorphism $g:\cW\ra \cW{}'$ representing
the germ $g:\Spec(\O_{\cZ{}',\bar{z}'}) \ra\Spec(\O_{\cZ,\bar{z}}) $,
and replace $\cU$ by some smaller \'etale neighborhood so that there is a
commutative diagram
$$
\begin{CD}
\cU @>>> \cW_0
@>>>W \\ 
@V\id VV @VVg_0 V @VVg V\\
\cU @>>> \cW{}'_0
@>>>W'.
\end{CD}
$$
Now let $\cU\subset \cU_1$ and $\cU\subset \cU_1{}'$ be the corresponding first order
extensions defined by $\cW$ and $\cW{}'$, respectively.
According to \cite{EGA IVd}, Theorem 18.1.2, we have an isomorphism
$\cU_1\simeq \cW_1\times_{\cW_1{}'} \cU_1{}'$, and conclude that the first order
extensions
$\cU_1,\cU_1{}'$ generate the same cyclic $\O_{\cX,\bar{x}}$-submodule in
$\shExt^1(\Omega^1_{\cX/k},\O_\cX)_{\bar{x}}$.
\qed

\medskip
Next we ask whether the collection of cyclic submodules in
$\shExt^1(\Omega^1_{\cX/k},\O_\cX)_{\bar{x}}$, $x\in |\cX|$ generated by
gtc-charts are the stalks of a coherent subsheaf. This is indeed true at
least over the subspace of nonnormality:

\begin{theorem}
\mylabel{injection}
Suppose $\cX$ admits a gtc-atlas $\shG$ with respect to $\sqM_\cX$ and $\rho$.
Let $\cD\subset \cX$ be the subspace of nonnormality, and $\shN_\cD\in\Pic(\cD)$ the
restricted conormal sheaf of $\shG$.
Then there is an injection $\shN_\cD^\vee\subset
\shExt^1(\Omega^1_{\cX/k},\O_\cX)_\cD$ whose stalks are the cyclic $\O_\cX$-submodules
generated by gtc-charts.
\end{theorem}

\proof
Choose gtc-charts $\cX\la \cU_i\ra \cZ_{i0}$ so that the disjoint union
$\cU=\bigcup \cU_i$ is an
\'etale covering of
$\cX$. Let $\cZ_i=\Spec(k[P_i])$ with $P_i=\sigma_i^\vee\cap M_i$ be the
corresponding Gorenstein toric variety, and $\rho_i\in P_i$ the element
defined by the Cartier divisor
$\cZ_{i0}\subset \cZ_i$. Then $\chi^{2\rho_i}\in k[P_i]$ defines a first order
extension
$\cZ_{i0}\subset \cZ_{i1}$, and by \cite{EGA IVd}, Theorem 18.1.2 
there is a cartesian diagram
$$
\begin{CD}
\cU_i @>>> \cU_{i1}\\
@VVV @VVV\\
\cZ_{i0} @>>> \cZ_{i1},
\end{CD}
$$
whose vertical arrows are \'etale.
We have to understand how the first order extensions
$\cU_i\subset \cU_{i1}$ differ on the overlaps $\cU_{ij}=\cU_i\times_\cX \cU_j$.
Fix a point $u\in \cU_{ij}$, and choose
an affine \'etale neighborhood $\cV{}'\ra \cU_{ij}$ of $u$
so that there is an isomorphism 
$\phi:U_j|_{\cV{}'}\ra U_i|_{\cV{}'}$ of
log spaces. Such isomorphism is given by a bijection
$\phi:p_1^*(\shM_{U_j})\ra p_0^*(\shM_{U_j})$.
Here $p_0$ and $p_1$ are the projections from $V'$ onto the
second and first factor of $\cU_{ij}$, respectively (compare Section 
\ref{cech cohomology}). Note that $\cV{}'=\cV{}'_{ijx}$ depends on $i,j,x$,
but we suppress this dependence to keep notations simple.

The sections $\rho_i\in\Gamma(\cU_i,\shM_{U_i})$ are lifts for
$\rho|_{\cU_i}\in\Gamma(\cU_i,\sqM_X)$, hence
$\phi(\rho_j|_{\cV{}'})= e'\cdot \rho_i|_{\cV{}'} $ for some
$e'\in\Gamma(\cV{}',\O_\cX^\times)$. Recall from Section
\ref{restricted conormal sheaf} that the restricted conormal sheaf $\shN_\cD$
is defined in terms of such $e'$.

Let $x\in |\cX|$, $z_i\in \cZ_i$, and $z_j\in \cZ_j$ be the images of 
$u\in \cU_{ij}$. In the proof of Proposition \ref{extension class}, we
constructed a bijection
$g:\Spec(\O_{\cZ_j,\bar{z}_j})\ra\Spec(\O_{\cZ_i,\bar{z}_i})$ inducing
a commutative diagram
$$
\begin{CD}
\Spec(\O_{\cX,\bar{x}}) @>>>\Spec(\O_{\cZ_{j0},\bar{z}_j})
@>>>\Spec(\O_{\cZ_j,\bar{z}_j}) \\ 
@V\id VV @VVg_0 V @VVg V\\
\Spec(\O_{\cX,\bar{x}}) @>>>\Spec(\O_{\cZ_{i0},\bar{z}_i})
@>>>\Spec(\O_{\cZ_{i},\bar{z}_i}).
\end{CD}
$$
By its very definition, the map $g$ sends 
$\chi^{\rho_j}\in\O_{\cZ_{j},\bar{z}_j}$ to 
$e'\cdot \chi^{\rho_i}\in\O_{\cZ_i,\bar{z}_i}$.
If follows that the extension class 
$\lambda_j\in\Gamma(\cU_j,\shExt^1(\Omega^1_{\cX/k},\O_\cX))$ of
$\cU_j\subset \cU_{j1}$ and the extension class
$\lambda_i\in\Gamma(\cU_j,\shExt^1(\Omega^1_{\cX/k},\O_\cX))$ of 
$\cU_i\subset U_{i1}$ are related
by $\lambda_j|_{\cV{}'}=e'\cdot\lambda_i|_{\cV{}'}$, at least after
refining $\cV{}'$.
This explains why the local extension classes $\lambda_i$ do not necessarily
satisfy the cocycle condition. However, we showed in Section \ref{restricted
conormal sheaf} below Equation (\ref{comparison}) that the cocycle condition
for $e'$ holds after restricting to the space of nonnormality $\cD\subset \cX$.

To be precise, set $\cV=\bigcup \cV{}'$, where the disjoint union runs over all
\'etale neighborhoods $\cV{}'=\cV{}_{i,j,x}'$. Then the
canonical map
$\cV\ra \cU\times_\cX \cU$ is an
\'etale covering. In this set-up, $\lambda\in
C^1(\cV/\cU,\shExt^1(\Omega^1_{\cX/k},\O_\cX))$ is a cochain.
Restricting to $\cD$ we obtain another cochain 
$\lambda_\cD\in C^1(\cV/\cU,\shExt^1(\Omega^1_{\cX/k},\O_\cX)\otimes\O_\cD)$.
On each $\cU_i$, the section $\rho_i\in\Gamma(\cU_i,\shM_{U_i})$ defines a
trivialization of
$\shN_\cD$, so we get an identification 
$$
C^1(\cV/\cU,\shExt^1(\Omega^1_{\cX/k},\O_\cX)\otimes\O_\cD)=
C^1(\cV/\cU,\shExt^1(\Omega^1_{\cX/k},\O_\cX)\otimes\shN_\cD).
$$
Now $\lambda_\cD$, viewed as a cochain with values in
$\shExt^1(\Omega^1_{\cX/k},\O_\cX)\otimes\shN_\cD)$, satisfies the cocycle
condition, according to the arguments below Equation (\ref{comparison}).
Consequently, $\lambda_\cD$ defines a global section of
$\shExt^1(\Omega^1_{\cX/k},\O_\cX)\otimes\shN_\cD$, and in turn the desired
homomorphism
$\shN_\cD^\vee\ra \shExt^1(\Omega^1_{\cX/k},\O_\cX)_\cD$. 
A local computation shows that this map is
bijective in codimension $\leq 1$. Here we use the assumption that our toric
varieties
$Z_i$ are regular in codimension $\leq 2$.
Since $\cD$ has no embedded component by Ishida's Criterion
(\cite{Oda 1988}, page 126), the map 
$\shN_\cD^\vee\ra\shExt^1(\Omega^1_{\cX/k},\O_\cX)$ is injective everywhere.
\qed

\medskip
For gtc-atlases, the restricted conormal sheaf $\shN_\cD$ thus has
two interpretations. First in terms of cocycles obtained
from $\rho\in\Gamma(\cX,\sqM_\cX)$ as in Section
\ref{restricted conormal sheaf}, and second in term
of first order extensions $\cU\subset \cU_1$ as in Theorem \ref{injection}.
We now state a generalization of Kato's result, who considered
spaces with normal crossing singularities (\cite{Kato 1996}, Theorem 11.7):

\begin{theorem}
\mylabel{gtc existence}
Let $\shG$ be a gtc-atlas on $\cX$ with respect to $\sqM_\cX$ and
$\rho\in\Gamma(\cX,\sqM_\cX)$, and $\cD\subset \cX$ the space of
nonnormality. Then there is a global log space $X\in\shG$ if and
only if the restricted normal sheaf
$\shN^\vee_\cD\subset\shExt^1(\Omega^1_{\cX/k},\O_\cX)$ extends to an invertible
$\O_\cX$-module.
\end{theorem}

\proof
This is Theorem \ref{existence} in the special case of gtc-atlases.
\qed

\begin{remark}\label{existence of gerbe}
In the normal crossings case the gerbe $\shG$ on $\cX$ is uniquely
determined by the requirement that $\overline\shM_\cX = f_*\NN_\cS$
for $f:\cS\to \cX$ the normalization. This is due to the fact that
such log structures are locally unique as shown in \cite{Kato 1996},
see also \cite{Kawamata; Namikawa 1994}. Indeed, if $\cX\to \Spec
k[z_1,\ldots,z_n]/(z_1\ldots z_r)$ is \'etale in $x\in|\cX|$ then
there exist $m_1,\ldots,m_r\in \shM_{\cX,\bar x}$ generating
$\overline\shM_{\cX,\bar x}= \NN^r$ and with $\alpha_X(m_i)= z_i$,
$i=1,\ldots,r$. For any other choices $m'_1,\ldots,m'_r$ the map
$m_i\mapsto m'_i$ defines uniquely an automorphism of $\shM_{\cX,\bar
x}$ fixing $\O^\times_{\cX,\bar x}$.

This argument does not work if $\overline\shM_{\cX,\bar x}$ has
relations. For example, consider the quadruple point $\cX=\Spec
k[z_1,z_2,z_3,z_4]/(z_1 z_3,z_2 z_4)$ with $\overline\shM_{\cX}$ the
ghost sheaf induced by the embedding into the toric variety $\Spec
k[z_1,z_2,z_3,z_4]/(z_1 z_3-z_2 z_4)$, with $k$ separably closed.
Then the set of isomorphism classes of log structures on $\cX$ is
canonically $(\NN\setminus \{0\})^4\times k^\times$, as explained in
\cite{Gross; Siebert 2003}, Example~3.13.  Note this example is
normal crossings away from the distinguished closed point of
multiplicity $4$ and hence the non-uniqueness is concentrated at this
point.

For the general case $\cX=\Spec k[P]/(\chi^\rho)$ with a Gorenstein
toric monoid $P$, $\rho\in P$ the distinguished element, and
$\overline\shM_\cX$ the ghost sheaf induced by the embedding into
$\Spec k[P]$, Proposition~3.14 in the same paper says the following.
Let $\underline x\in\cX$ be the distinguished closed point. Then the
set of isomorphism classes of germs at $\underline x$ of
gtc-structures on $\cX$ with ghost sheaf $\overline\shM_\cX$ injects
into $\shExt^1 (\overline\shM^\gp_\cX, \O_\cX^\times)_{\underline x}$ by
associating the extension class. Moreover, there is an explicit
description of both $\shExt^1 (\overline\shM^\gp_\cX,
\O_\cX^\times)_{\underline x}$ and the image of the germs of log
structures in terms of functions $h_p$, $p\in P$, on open subsets of
$\cX$. The function $h_p$ is defined on the complement of
$V(\chi^p)\subset \cX$. Conversely, given $(h_p)_{p\in P}$ such that $h_p$
extends to $X$ by $0$ then $p\mapsto h_p$ defines a chart for the
corresponding log structure.

This description also suggests a notion of \emph{type} for germs of
log structures on $\cX$, namely if their representatives $(h_p)$,
$(h'_p)$ differ only by invertible functions \cite{Gross; Siebert
2003}, Definition~3.15. Globally two log structures are of the same
type if they are of the same type at each point. Log structures of
the same type have charts with image in the same toric variety and
inducing the same combinatorial identification of prime components
with toric prime divisors. In the example of the quadruple point
fixing the type means choosing an element in $(\NN\setminus\{0\})^4$.

Note that in any case $\overline\shM_\cX$ is naturally a subsheaf of
$f_*\NN_{\cS}$ for $f:\cS\to \cX$ the normalization, and this subsheaf
determines the type of log structure. Indeed, it suffices to check
this for $\cX$ the boundary divisor in a toric variety. Let $g_1:
\cX\to\Spec k[P]$, $g_2: \cX\to\Spec k[P]$ be isomorphisms of $\cX$ with
the boundary divisor of the toric variety $\Spec k[P]$ inducing the
same embedding of $P$ into $\NN^r$, $r$ the number of irreducible
components of $\cX$. Then for any $p\in P$ the orders of vanishing of
$g_1^*(\chi^p)$ and $g_2^*(\chi^p)$ along the toric prime divisors
agree and hence there exists $h_p\in \Gamma(\O_\cX)$ with
$g_1^*(\chi^p)=h_p g_2^*(\chi^p)$. This shows that the two log
structures induced by $g_1$ and $g_2$ are of the same type. In
particular, comparing the type of log structures for a given set of
charts is a finite problem that in practice can often be done by
hand.

Taken together this gives a three-step solution to the problem of
constructing gtc structures on a given algebraic space $\cX$: First
determine the type of gtc structure by covering $\cX$ with finitely
many charts of the same type on overlaps as discussed. In the next
step one needs to compare the selected sections of $\shExt^1
(\overline\shM_\cX^\gp, \O_\cX^\times)$ and adjust if necessary.
Although this step is still abelian in nature, it is probably the
most difficult one in practice. On the other hand, on the
(semi-stable) normal crossings locus where
$\overline\shM_\cX=f_*\NN_\cS$ the subsheaf of $\shExt^1
(\overline\shM_\cX^\gp, \O_\cX^\times)$ parametrizing log structures
of semi-stable type is trivial and hence has a unique section. This
follows from the mentioned explicit description of this sheaf in
\cite{Gross; Siebert 2003}, and it reflects the uniqueness of the
gerbe $\shG$ on such spaces discussed above. Thus this second step is
simple on the normal crossings locus. The third and last step is an
application of the theorem above.

That this is indeed a viable approach has been shown in \cite{Gross;
Siebert 2003}. In this paper $\cX$ is a union of toric varieties and
the result is a classification of gtc structures in terms of a
certain, computable sheaf cohomology group on a real integral affine
manifold $B$ built on the dual intersection complex of $\cX$. In this
case the given cell decomposition of $B$ already determines the ghost
sheaf.
\end{remark}

\section{Triple points and quadruple points}

It is now time to illustrate the general theory with some concrete
examples. The examples are normal crossing except at finitely many
points. According to Remark~\ref{existence of gerbe} to define the
gerbe $\shG$ it suffices to specify charts at these points.

\begin{example}
We start by looking at 3-dimensional affine toric varieties
$Z=\Spec k[\sigma^\vee\cap \ZZ^3]$ that are Gorenstein and $(R_2)$, such
that the boundary divisor $Z_0$ has three irreducible components.
Let $\rho\in \sigma^\vee\cap \ZZ^3$ be the unique element
with $\rho+\sigma^\vee\cap \ZZ^3=(\interiour\sigma^\vee)\cap\ZZ^3$.
After changing coordinates, we may assume that
$\rho=(0,0,1)$. Let $H\subset \RR^3$ be the affine hyperplane
defined by the affine equation $\rho^\vee=1$.
Then the cone $\sigma$ is generated by a \emph{lattice triangle} in $H$
generated by $v_1,v_2,v_3\in H$ such that the vertices are the only boundary
lattice points.

Applying an integral linear coordinate change fixing $\rho\in\ZZ^3$,
we may assume $v_1=(0,0,1)$, $v_2=(1,0,1)$
and $v_3=(a,b,1)$ for some $a,b\in\ZZ$.
Making further coordinate changes using the matrices
$$
\begin{pmatrix}
1\\
&\pm 1\\
&& 1\\
\end{pmatrix}
\quadand
\begin{pmatrix}
1&\pm 1\\
&1\\
&&1\\
\end{pmatrix}
$$
we end up with $0\leq b$ and $0\leq a<b$. The condition that the
segments $\overline{v_1v_3}$ and $\overline{v_2v_3}$ contain no
additional lattice point means that both $a,a-1$ are prime to $b$.
Note that $b$ is necessarily odd, because either $a$ or $a-1$ is
even. Moreover, $b\geq 3$ implies $a\geq 2$. The case $v_3=(0,1,1)$
yields the regular toric variety. The simplest nontrivial case is
therefore $v_3=(2,3,1)$, which defines the unique isomorphism class
of lattice triangle with one interior lattice point and three
boundary lattice points.

The boundary divisor $Z_0\subset Z$ decomposes into three
irreducible  components $Z_0=Z_{01}\cup Z_{02}\cup Z_{03}$
corresponding to the  vectors $v_1,v_2,v_3$. Each $Z_{0i}$ is a
2-dimensional affine toric variety. Its cone is the image of $\sigma$
under the canonical projection $\ZZ^3\ra\ZZ^3/\ZZ v_i$. Since
$\det(v_1,v_2,v_3)=b$, the $Z_{0i}$ are affine toric surfaces
containing the  rational Gorenstein singularity of type $A_{b-1}$.
Note that the underlying scheme $\cZ_0$ is determined up to
isomorphism by the integer $b\geq 1$. This is because the
normalization map $\amalg \cZ_{0i}\ra \cZ_0$ is determined in
codimension $\leq 1$, compare the discussion in \cite{Reid 1994},
Section 2.

On the other hand, the log space $Z_0$ depends on the integer $a$.
How many such $a$ are possible? Suppose for a moment that $b=p^n$ is
an odd prime power. Then both $a,a-1$ are prime to $p$ if and only if
$a$ is neither in $p\ZZ/(p^n)$ nor in $1+p\ZZ/(p^n)$. Hence there are
$p^n-2p^{n-1}=p^{n-1}(p-2)$ choices for $a$. In general, decompose
$b=\prod p_i^{n_i}$ into prime factors. Then there are $\prod
p_i^{n_i-1}(p_i-2)$ possibilities for $a$.

Now suppose we have a 2-dimensional algebraic $k$-scheme $\cX$ that
is normal crossing in codimension $\leq 1$ and whose irreducible
components $\cX_i$ are normal. Let $x_j\in \cX$ be the
closed points where at least three irreducible components meet. Away
from the $x_j$ our gtc-atlas is uniquely determined by the
requirement that $\overline\shM_\cX$ agree with $f_*\NN_\cX$,
$f:\cS\to\cX$ the normalization. We assume that each closed point
$\underline x_j\in \cX$ that is not normal crossing is \'etale
locally isomorphic to $\cZ_0=\cZ_{j0}$ at the origin for certain odd
integers $b_j\geq 1$. The choice of integers $0\leq a_j<b_j$ such
that both $a_j,a_j-1$ are prime to $b_j$ now specifies a gtc-chart
$\shG$ on $\cX$ that is naturally compatible with the already chosen
gtc-atlas on the complement of the $\underline x_j$.
\end{example}

\begin{example}
Let us now consider another example.
Let $Z=\Spec k[\sigma^\vee\cap \ZZ^3]$ be a 3-dimensional Gorenstein toric
variety satisfying $(R_2)$, such that the boundary divisor $Z_0$
has four irreducible components.
Now the cone $\sigma\subset N\otimes\RR$ is generated by a \emph{lattice
tetragon} in the affine hyperplane $H\subset N\otimes\RR$
whose vertices are the only boundary lattice points.
Let 
$v_1,\ldots,v_4\in H$ be the vertices of such a lattice tetragon. After an
integral coordinate change, we may assume
$v_1=(0,0,1)$, $v_2=(1,0,1)$, $v_3=(a,b,1)$ with $0\leq a<b$ and
$\gcd(a-1,b)=1$, and $v_4=(c,d,1)$ with $\gcd(c,d)=\gcd(c-a,d-b)=1$.
The convexity condition is $ad-bc>0$ and $d>0$.

Let $Z_{01},\ldots,Z_{04}\subset Z_0$ be the irreducible components
corresponding to the vectors
$v_1,\ldots,v_4\in\sigma$, respectively. Each $Z_{i0}$ is a Gorenstein toric
variety. We have $\det(v_4,v_1,v_2)=d$, so the invariant closed point on
$Z_{01}$ is the rational Gorenstein singularity of type $A_{d-1}$.
Similarly, $Z_{02}$ has type $A_{b-1}$, and $Z_{03}$ has type
$A_{b-d+ad-bc-1}$, and
$Z_{04}$ has type
$A_{ad-bc-1}$.

Let us now concentrate on the special case $a=b=d=1$ and $c=0$, that
is $v_3=(1,1,1)$ and  $v_4=(0,1,1)$. This corresponds to the unique
lattice tetragon containing precisely four lattice points. Then every
irreducible component $Z_{i0}$ is smooth. The boundary divisor $Z_0$
is a complete intersection isomorphic to the spectrum of
$A=k[x,y,u,v]/(xy, uv)$.  Note that we may view $Z_0$ as the product
of two 1-dimensional normal crossings. The space of nonnormality
$D\subset Z_0$ is the union of the four coordinate axis in $\AA^4_k$,
given by the subring in $k[x]\times k[y]\times k[u]\times k[v]$ of
polynomials with identical constant term. Using the coordinates
$x,y,u,v$, we calculate
\begin{align*}
\Ext^1(\Omega^1_{Z_0/k},\O_{Z_0})&= 
A/(\frac{\partial}{\partial x} xy,\frac{\partial}{\partial y} xy)
\oplus
A/(\frac{\partial}{\partial u} uv,\frac{\partial}{\partial v} uv)\\
&=k[u,v]/(uv)\oplus k[x,y]/(xy).
\end{align*}
Under this identification, the restricted conormal sheaf
$\shN_D\subset\shExt^1(\Omega^1_{Z_0/k},\O_{Z_0})$ corresponds to the
diagonal submodule $(f(u,v),f(x,y))$.

Here is an example for a proper algebraic surface having such a
quadruple point: Let $S\ra\PP^1_k$  be a Hirzebruch surface of degree
$e\geq 0$. We denote by $C_1$ the unique section with $C_1^2=-e$, and
choose another section $C_3\subset S$ with $C_2^2=e$. Let
$C_2,C_3\subset S$ be the fibers over $0,\infty\in\PP^1_k$,
respectively. Then $C=C_1\cup C_2\cup C_3\cup C_4$  forms a 4-cycle
of smooth rational curves. Now choose an isomorphism $C_2\ra C_4$ 
sending $C_1\cap C_2,C_2\cap C_4$ to $C_1\cap C_4,C_4\cap C_3$,
respectively, and let $C_1\ra C_3$ be a similar isomorphism. Then
define $X$ to be the proper algebraic space obtained from $S$ by
identifying $C_1,C_3$ and $C_2,C_4$ with respect to these maps. Then
$X$ has normal crossing singularities except for a single closed
point $x\in |\cX|$, whose preimage on $S$ are the nodal points of
$C$. \'Etale locally near $x$, the space $X$ is isomorphic to the
boundary divisor $Z_0$. Hence, as in the prevous example, $X$ is
endowed with a gtc-atlas $\shG$ with the property that the ghost
sheaf agrees with $f_*\NN_\cS$ away from $x$, $f:\cS\to\cX$ the
normalization. We examined similar surfaces in connection with
degenerations of primary Kodaira surfaces \cite{Schroeer; Siebert
2002}.
\end{example}

\section{Smooth log atlases}

In this short section we propose a tentative generalization of
gtc-atlases using the concept of smoothness in the category
of log spaces.
Recall that a morphism $f:X\ra Y$ of fine log spaces is called
\emph{smooth} if, \'etale locally, there are charts $P_X\ra\shM_X$,
$Q_Y\ra\shM_Y$, and
$Q\ra P$ for $f$ such that the induced morphism 
$\cX\ra \cY\otimes_{\ZZ[Q]}\ZZ[P]$ of algebraic spaces is \'etale,
and that kernel and the torsion part of the cokernel for $Q^\gp\ra P^\gp$
are groups of order prime to the characteristic of the ground field.
Equivalently, the morphism $f:X\ra Y$ satisfies the lifting
criterion for log Artin rings similar to the classical lifting criterion for
smoothness of schemes. 
It turns out that smooth log spaces behave very much like smooth
spaces, and can be treated with similar methods.
For more details on smooth morphism of log spaces
we refer to 
\cite{Kato 1989}, Section 3. 

We now consider the following situation. 
Fix a ground field $k$ and a fine monoid $Q$.
Let $(\Spec(k),Q)$ be the log structure associated to the prelog structure
$$
Q\lra k,\quad q\mapsto\begin{cases}
1\quad\text{if $q=0$,}\\
0\quad\text{otherwise}.
\end{cases}
$$
The geometric stalk of $\shM_{(\Spec(k),Q)}$ is $(k^\sep)^\times\oplus Q$.
Now let $\cX$ be an algebraic $k$-space of finite type endowed
with a constructible monoid sheaf $\sqM_\cX$ with fine stalks.
We also assume that we have a fixed monoid homomorphism $\rho:Q\ra\sqM_\cX$.
We propose the following definition:

\begin{definition}
\mylabel{smooth atlas}
A log atlas $\shG$ on $\cX$ with respect to $\sqM_\cX$ is called \emph{smooth}
if there is an \'etale covering $\cU\ra \cX$,
a log space $U\in\shG$, and a smooth morphism of log spaces
$U\ra (\Spec(k),Q)$ compatible with $\rho:Q\ra\sqM_\cX$.
\end{definition}

Note that a morphism $U\ra(\Spec(k),Q)$ compatible with $\rho$ is
nothing but a lifting $\tilde{\rho}:Q\ra\shM_U$ of $\rho:Q\ra\sqM_U$, thanks
to the splitting of
$\shM_{(\Spec(k),Q)}$.
Observe that gtc-atlases are smooth log atlases: In this special case we
have $Q=\NN$, and the fixed morphism $\rho:Q\ra\sqM_\cX$ corresponds to the
fixed section
$\rho\in\Gamma(\cX,\sqM_\cX)$.

We expect that the notion of smooth log atlases will be 
crucial in studying degenerations and deformations over
higher dimensional base schemes.

\section{Kato fans}

In this section we recall a combinatorial object introduced by Kato
\cite{Kato 1994} under the name \emph{fan}. To avoid confusion with
toric geometry, we shall use the  term \emph{Kato fan}. This concept
will be a convenient framework for our mirror construction in the
next two sections. To keep the discussion within limits we work in
the category of schemes locally of finite type over a ground field
$k$ rather than algebraic spaces. Also the log structures are now
defined on the Zariski site. See \cite{Niziol 2003} for a detailed
comparison between log structures on the Zariski and \'etale sites.
Essentially this only rules out self-intersecting components in our
construction, confer Kato's comment in \cite{Kato 1994}, Remark~1.8. 
We may avoid this restriction with a little more effort, confer
\cite{Gross; Siebert 2003}, Section~2.2.

Recall from \cite{Kato 1994}, Definitions~9.1 and 9.3 that a monoidal
space is a topological space $T$ endowed with a sheaf of sharp
monoids  $M_T$, and that a \emph{Kato fan} is a monoidal space
$(T,M_T)$  that is locally of the form
$$
(\Spec(P),M_{\Spec (P)}),
$$
where $\Spec(P)$ is the set of prime ideals in some monoid $P$. 
Here the notation is adopted from commutative algebra.
In multiplicative notation, $I\subset P$ is an ideal if
$PI\subset I$, and it is a prime ideal if $P\setminus I$ is a submonoid of
$P$ (\cite{Kato 1994}, Definition~5.1). The spectrum
$\Spec(P)$ is the set of prime ideals in $P$ with the topology
generated by $D(f)=\{\mathfrak{p} \in\Spec(P)\,|\, f\not\in
\mathfrak{p}\}$ for $f\in P$. The sections of $M_{\Spec (P)}$ over
$D(f)$ are
$$
S^{-1}P/(S^{-1}P)^\times\quad \text{with}\ S=\{f^n\,|\, n\ge0\}.
$$
Similarly, for a prime ideal $\mathfrak{p}\subset P$ we write
$
P_\mathfrak{p}= S^{-1}P/(S^{-1}P)^\times,
$
with $S=P\setminus \mathfrak{p}$. This is the stalk of $M_{\Spec(P)}$
at $\mathfrak{p}$.

The affine Kato fan $\Spec(P)$ is finite if $P$ is finitely
generated. A Kato fan $T$ is \emph{locally of finite type} if the
monoids $P$ can be chosen to be finitely generated. In contrast to
the situation in \cite{Kato 1994} we will not be able to restrict to
integral monoids as we will see shortly. A Kato fan that is locally
of finite type is locally finite. A convenient way to think about 
locally finite topological spaces is as partially ordered sets via
$$
x\leq y\quad\Longleftrightarrow\quad x\in \overline{\left\{y\right\}}.
$$
Reversing this ordering leads to the \emph{dual space $F^*$}. In other words,
$F^*=F$ as sets, but $U\subset F^*$ is open iff $U \subset
F$ is closed. A sheaf $P$ on $F$ is equivalent to a set of
monoids $P_x$ indexed by $x\in F$, together with a compatible system
of generization maps
$\varphi_{yx}: P_x\ra P_y$
for any $x\leq y$. 

Kato fans arise in log geometry as follows. For a scheme $\cX$ with
fine log atlas $\shG$ and $x\in \cX$ denote by $I(\shG,x)\subset
\mathfrak{m}_x$ the ideal generated by the image of $P\setminus
\alpha_x^{-1}(\O_{\cX,x}^\times)$ for any chart $\alpha_x: P\to
\O_{\cX,x}$ at $x$. Note that $I(\shG,x)$ depends only on $\shG$ and
not on the particular chart. We are interested in equivalence classes
of log structures with the same ghost sheaf and the same set of ideals
$\shI(\shG,x)$.

\begin{definition}\label{pre-gtc atlas}
Let $\cX$ be a scheme endowed with a sheaf of fine sharp monoids
$\sqM_\cX$. Suppose we have an  \'etale covering $\cU_i\to \cX$ and log
spaces  $U_i$ together with an
identification $\sqM_{U_i}\simeq \sqM_{\cU_i}$.
Let $q_i:\shM_{U_i}\to \overline\shM_{\cU_i}$ be the quotient map.
We call  $(U_i,q_i)$ a \emph{pre-gtc atlas} if:
\begin{enumerate}
\item 
For each $i$ there exists an \'etale map $\cU_i\to \cZ=\Spec
k[\sigma^\vee\cap M]/(\chi^{\rho_i})$ to the reduced boundary divisor
of a Gorenstein toric variety inducing the log space $U_i$.
\item 
For any $x\in \cX$ and $p\in \overline\shM_{\cX,x}$ the ideal
$\shI(U_i,x)\subset \O_{\cX,x}$ generated by $\alpha_i( q_i^{-1}(p))$ is independent of
the choice of $i$ with $x\in \cU_i$.
\end{enumerate}
\end{definition}

In the situation of the definition the pull-backs of $\rho_i\in
\Gamma(\cZ,\sqM_\cZ)$ glue to a distinguished section
$\rho\in\Gamma(\cX, \sqM_\cX)$. This is true because for a Gorenstein
sharp toric monoid $P$ there is a unique element $\rho\in P$ with the
property $P\setminus (\rho+P)= \partial P$. Moreover, for any $x\in
\cX$ there is a well-defined ideal $\shI(\overline\shG,x)\subset
\mathfrak{m}_x$ by taking $\shI(U_i,x)$ for any $i$ with $x\in
\cU_i$. A scheme with a pre-gtc atlas induces a Kato fan (cf.\
\cite{Kato 1994}, Proposition~10.1 for an analogue for toroidal
varieties):

\begin{proposition}
\label{associated fan}
Let $\cX$ be a scheme endowed with a sheaf of fine sharp monoids
$\sqM_\cX$, together with a pre-gtc atlas $\overline\shG$. Let
$\rho\in\Gamma(\cX,\sqM_\cX)$ be the distinguished section. Then
\begin{enumerate}
\item 
The ideal
$I(\overline\shG,x)\subset \O_{\cX,x}$ is a prime ideal for every $x\in \cX$.
\item
The set
$F(\cX)= \{x\in \cX\mid I(\overline\shG,x)=\mathfrak{m}_x\}$ endowed with
the subspace topology from $\cX$ and the monoid sheaf
$M_{F(\cX)}=\overline{\shM}_\cX/(\rho)|_{F(\cX)}$ is a Kato fan
locally of finite type.
\item
There is a morphism $ \pi:(\cX,\sqM_\cX \dd (\rho))\ra
(F(\cX),M_{F(\cX)})$ mapping $x\in \cX$ to the point of
$F(\cX)\subset \cX$ corresponding to the prime ideal
$I(\overline\shG,x)\subset\O_{\cX,x}$, and the canonical map
$\pi^{-1}M_{F(\cX)}\ra\overline{\shM}_\cX\dd(\rho)$ is bijective.
\end{enumerate}
\end{proposition}

\begin{remark}
\label{monoid quotients}
In the statement of the proposition we are taking a certain quotient of a
monoid $M$ by an ideal $J=(\rho)=\rho+ M$. This quotient is defined
as the set consisting of $M\setminus J$ together with one more point
$\infty$. For $m,m'\neq0$ set $m+ m'=\infty$ if one of $m,m'$ equals
$\infty$, or if $m+ m'\in J$ as sum in $M$. Otherwise the sum $m+ m'$
is taken in $M$.

This construction has the following categorical meaning. Consider the
category of monoid homomorphisms $M\to M'$ mapping $J$ to an
attractive point $\infty\in M'$, that is, with $\infty+m=\infty$ for
all $m\in M'\setminus \{0\}$. Then
\[
\varphi: M\lra  M\dd J,\quad
m\longmapsto \left\{
\begin{array}{ll}
m,& m\in M\setminus J\\
\infty,&m\in J,
\end{array}\right.
\] 
is an initial object in this category. Note that unless $J=0$ our
quotients are never integral and that any ideal in $M\dd J$ contains
$\infty$. Such ideal quotients in the category of monoids are
compatible with ideal quotients in the category of rings in the
following sense. Let $(\chi^J)\subset k[M]$ be the ideal generated by
monomials $\chi^m$ with $m\in J$. Then there is a canonical
isomorphism $k[M\dd J]/(\chi^\infty)= k[M]/(\chi^J)$.

The referee pointed out that this is indeed not a proper quotient in the
category of monoids, which is why we use the double slash notation.
\end{remark}

\noindent
\emph{Proof of Proposition \ref{associated fan}.}
Because the problem is local we may restrict ourselves to the case
that $\cX$ has an \'etale morphism to
$\Spec k[P]/(\chi^\rho) = \Spec k[P\dd (\rho)]/(\chi^\infty)$ 
inducing the log structure, where
$P=\sqM_{\cX,x}$. For $\Spec k[P]$, Kato proved
 statements (1)--(3) in \cite{Kato 1994}, Section~10. In particular, for
each prime ideal $\mathfrak{p}\subset P$ there is exactly one point
$x\in\Spec k[P]$ such that $P\setminus \mathfrak{p}$ generates the
maximal ideal at $x$, and conversely. Therefore the points $x\in
F(\cX)$ are in one-to-one correspondence with prime ideals in $P$
contained in $P\setminus (\rho+P)$. Hence $F(\cX) =\Spec (P\dd (\rho))$.
Statement~(3) follows from the corresponding statement for
$\Spec k[P]$ by dividing out the ideal $(\rho)$.
\qed
\medskip

The Kato fan $(F(\cX),M_{F(\cX)})$ in Proposition \ref{associated fan} is
a hull for $(\cX, \sqM_\cX\dd (\rho))$ rather than for $(\cX,
\sqM_\cX)$. For our construction in the next section we need
an additional structure on $F(\cX)$ coming from the sheaf
$\sqM_\cX$ on $\cX$.

\begin{definition}
A \emph{gtc-structure} on a monoidal space $(F,M_F)$ is a sheaf $P$
of Gorenstein sharp toric monoids, together with an isomorphism $P\dd
(\rho+P) \simeq M_F$ for $\rho\in\Gamma(F,P)$ the distinguished
section. A \emph{gtc-fan} is a Kato fan with a gtc-structure. The
notation will be $(F,P,\rho)$. \end{definition}

If $P$ is a Gorenstein sharp toric monoid with distinguished element
$\rho$ then the restriction of $M_{\Spec P}$ to $\Spec(P\dd (\rho))$ is
a gtc-structure on $(\Spec(P\dd (\rho)), M_{\Spec(P\dd (\rho))})$. Hence
the following is a direct consequence from the proof of 
Proposition~\ref{associated fan}.

\begin{proposition}
The Kato fan $(F(\cX), M_{F(\cX)})$ from Proposition~\ref{associated fan}
has a gtc-structure $(P_{F(\cX)}, \rho)$.
\end{proposition}

We call $(F(\cX),P_{F(\cX)},\rho)$ the \emph{gtc-fan associated to
$(\cX,\overline\shM_\cX,\rho)$}.  Next we show how to construct a space
with toric components from a gtc Kato fan. For a toric monoid
$P=\sigma\cap \ZZ^d$, $P^\gp=\ZZ^d$, there is a one-to-one
correspondence between faces $\tau$ of $\sigma$ and those submonoids
$Q\subset P$ whose complement is a prime ideal, by taking the
integral points of $\tau$. Such submonoids are commonly called
\emph{faces} of $P$. Its (co-) dimension is the (co-) dimension of
$\tau$ in $\sigma$. Faces of codimension $1$ are \emph{facets}. We
write $P^\vee$ for the dual monoid $\Hom(P,\NN)$.

Let $(F,P,\rho)$ be a gtc Kato fan. For any $x\in F$ we have
the ring $k[P_x^\vee]$. Evaluation at $\rho\in \Gamma(P)$
defines a grading $P_x^\vee\to \NN$. We thus obtain a projective
scheme
$$
\cY_x=\Proj ( k[P_x^\vee]).
$$
The generization maps for the stalks of $P$ tell how to glue these
spaces according to the following lemma.

\begin{lemma}\label{closed embeddings}
For any toric monoid $P$ and $\mathfrak{p}\in \Spec(P)$ there exists
a canonical surjective morphism
$k[P^\vee]\ra k[P_\mathfrak{p}^\vee]$.
These morphisms are natural with respect to inclusion of
prime ideals.
\end{lemma}

\proof
Let $S= P\setminus \mathfrak{p}$ be the face associated to
$\mathfrak{p}$. As the elements of $S$ are invertible in $S^{-1}P$
the homomorphism
$P\ra P_\mathfrak{p}$
is surjective. Dualizing gives an injection
$
P_\mathfrak{p}^\vee\ra P^\vee$.
The image comprises those $\varphi:P\to \NN$ with $\varphi(S)=0$,
because $S^\gp$ is the kernel of $P^\gp\to P^\gp_\mathfrak{p}$.
Therefore $P^\vee\setminus P_\mathfrak{p}^\vee$ is an ideal. Letting
$I\subset k[P^\vee]$ be the ring-theoretic ideal generated by $\chi^m$
with $m\in P^\vee\setminus P_\mathfrak{p}^\vee$, we obtain the desired
surjection
$$
k[P^\vee]\lra k[P^\vee]/I= k[P_\mathfrak{p}^\vee].
$$
If $\mathfrak{p}\subset\mathfrak{q}$ there is a factorization
$
\varphi_\mathfrak{q}: k[P^\vee]
\ra k[P^\vee_\mathfrak{p}]
\ra k[P^\vee_\mathfrak{q}]$, and this  gives naturality.
\qed

\medskip
For $x\leq y$ there exists a prime ideal $\mathfrak{p}$ of $P=P_x$
and an isomorphism $P_y \simeq P_\mathfrak{p}$ such that
$\varphi_{yx}: P_x\to P_y$ is the localization map $P\to
P_\mathfrak{p}$. This follows because locally around $y$ the monoidal
space $(F,P)$ is isomorphic to
$$
(\Spec (P_y\dd (\rho)), M_{\Spec (P_y)}|_{\Spec (P_y\dd (\rho))}).
$$
So we can apply Lemma \ref{closed embeddings}. The
epimorphism $q_{yx}: k[P_x^\vee] \to k[P_y^\vee]$ thus obtained
respects the grading. For any $x\leq y$ we therefore get a closed
embedding
$\varphi_{xy}: \cY_y\ra \cY_x$.
By compatibility with localization the $\cY_x$, $x\in F$, together with
the closed embeddings $\varphi_{xy}$ form a directed system of
projective toric schemes.

\begin{lemma}
\label{inverse limit}
The direct limit $\dirlim \cY_x$ exists as a reduced $k$-scheme locally
of finite type, and the maps $\cY_x\to \dirlim \cY_x$ are closed
embeddings. If $F$ is finite then $\dirlim \cY_x$ is projective.
\end{lemma}

\proof
We may assume that $F$ is finite. If there
is only one closed point $z\in F$, the direct limit is $\cY_z$, because
the $\varphi_{xy}$ are closed embeddings. In the general case, fix a
closed point $z\in F$, let $F_1\subset F$ be the set of points that
are generizations of $z$, and let $F_2\subset F$ be the set of points
that are generizations of a closed point different from $z$. Let
$\cY_1,\cY_2,\cY_{12}$ be the direct limits corresponding to
$F_1,F_2,F_1\cap F_2$, respectively. These are projective schemes by
induction on the cardinality of $F$. We now view $\cY=\dirlim \cY_x$ as a
coproduct $\cY_1\amalg_{\cY_{12}} \cY_2$.  According to \cite{Artin 1970},
Theorem 6.1, the coproduct exists as a reduced algebraic space over
$k$, with $\cY_i\to \cY_1\amalg_{\cY_{12}} \cY_2$ closed embeddings with
images covering $\cY_1\amalg_{\cY_{12}} \cY_2$ set-theoretically. Repeating
this construction with the compatible system of ample line bundles
$L_x\ra \cY_x$ corresponding to the ample invertible sheaves
$\O_{\cY_x}(1)$, we infer that the algebraic space $\cY$ carries a line
bundle whose restriction to each irreducible component is ample.
Hence $\cY$ is a projective scheme.
\qed
\medskip

We write $\cY_{(F,P,\rho)}=\dirlim \cY_x$.

\section{A naive mirror construction}

Let $(F,P,\rho)$ be a gtc Kato fan, and $\cY_{(F,P,\rho)}=\dirlim \cY_x$
the corresponding projective scheme from Lemma~\ref{inverse limit}. 
Our next goal is to define a pre-gtc atlas on $\cY_{(F,P,\rho)}$.
This requires some additional data leading to a selfdual structure,
which in turn gives a baby version of mirror symmetry.

First note that we have a canonical identification
$\Spec(P)^*\simeq\Spec(P^\vee)$ for any toric monoid
$P=\sigma\cap\ZZ^d$, by sending $\tau\cap P$ to $(\RR\tau)^\perp\cap
P^\vee$. We exploit this as follows: For any closed point $x\in F$
there is a continuous map
$$
\cY_x=\Spec( k[P_x^\vee])\lra
\Spec(P_x^\vee)\simeq(\Spec(P_x))^*\subset F^*.
$$
The collection of these maps descends to a continuous map
$\cY_{(F,P,\rho)}\ra F^*$. This map should come from a pre-gtc
atlas on $F^*$. Thus one ingredient to define the desired pre-gtc
atlas on $Y_{(F,P,\rho)}$ will be a monoid sheaf $Q$ over $F^*$ with
section $\rho^*\in \Gamma(F^*,Q)$ making $(F^*,Q,\rho)$ into a gtc Kato
fan.  Of course, we also need a compatibility condition relating
$(F,P,\rho)$ to $(F^*,Q,\rho^*)$. We call a map  $\lambda:Q\to A$
from a monoid into an abelian group \emph{affine} if
$\lambda-\lambda(0)$ is a homomorphism of monoids. 

\begin{definition}\label{duality datum}
A \emph{gtc duality datum} consists of the following:
\begin{enumerate}
\item
A gtc Kato fan $(F,P,\rho)$.
\item
A gtc Kato fan $(F^*,Q,\rho^*)$.
\item
A \emph{compatibility datum} between $(F,P,\rho)$ and
$(F^*,Q,\rho^*)$ as follows: For any generic point $y\in F$
and any closed point $x\leq y$ we have an affine injection
$$
\lambda_{xy}: Q_y^\vee\lra P_x^\gp
$$
identifying $Q_y^\vee$ with the cone with vertex $\rho(x)$ over a
subset of the face of $P_x^\gp$ corresponding to $y\in\Spec(P_x)\subset F$,
such that the one-dimensional face of $Q_y^\vee$ containing
$\lambda_{xy}^{-1}(0)$ equals $Q_x^\vee\subset Q_y^\vee$.
\end{enumerate}
The notation will be $(F,P,Q,\lambda=\{\lambda_{xy}\})$.
\end{definition}

The following picture illustrates compatibility data. The
four long dash-dotted lines are the rays of $P_x$, so we are looking
from inside $P_x$. Let the facet of $P_x$ containing the polygon $\Delta$
correspond to $y\in F$. Then the indicated cone over $\Delta$
represents $Q_y^\vee$.
\\[3ex]
\centerline{\includegraphics{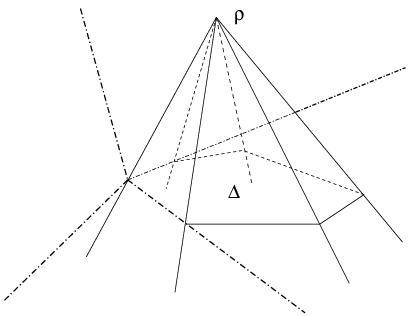}}\\[-2ex]
\mbox{}\hspace{2.5cm} Fig.\ 1\\[3ex]
Note that a compatibility datum is uniquely defined by the system of
polytopes $\Delta$, one for each facet of any $P_x$ with $x\in F$ a
closed point. The compatibility between the polytopes is best
expressed by saying that they give rise to a sheaf $Q$ on the dual
space $F^*$.
\medskip

We now explain how a gtc-duality datum gives rise to a pre-gtc
atlas on $\cY_{(F,P,\rho)}$.

\begin{construction}
(Construction of pre-gtc structure.)
We shall cover $\cY_{(F,P,\rho)}$ by divisors in affine toric
schemes, one for each generic point $y\in F$. Let $x\leq y$ be a
closed point. To $y\in\Spec(P_x)\subset F$ 
belongs a facet $S\subset
P_x$. Let $w\in P_x^\vee$ be the generator of the one-dimensional
face dual to $S$. Denote by $(P_x^\vee)_{(w)}$ the submonoid of
$(P_x^\vee)^\gp$ of terms of the form $p-a\cdot w$, $p\in P_x^\vee$,
$a\in\ZZ$ with $p(\rho)=a\cdot w(\rho)$. The notation comes from
interpreting $(P_x^\vee)_{(w)}$ as homogeneous localization of
$P_x^\vee$ with respect to the grading defined by 
$\rho$. The injection
$$
\lambda_{xy}-\rho:Q_y^\vee\lra P_x^\gp
$$
induces a bijection of groups $(Q_y)^\gp\simeq (P_x^\vee)^\gp$. We view
$Q_y$ as submonoid of $(P_x^\vee)^\gp$ via this bijection. With
this understood we have
\begin{eqnarray}\label{identification of cones}
(P_x^\vee)_{(w)}= Q_y\cap \rho^\perp.
\end{eqnarray}
Indeed, if $p-aw\in (P_x^\vee)_{(w)}$ then $(p-aw)(\rho)=0$ by
definition. To check that the image is in $Q_y$ it suffices to
evaluate its $\RR$-linear extension on $v-\rho$, for all vertices $v$
of the polygon $\Delta\subset S^\gp\otimes_\ZZ \RR$ spanning
$Q_y^\vee$:
$$
(p-aw)(v-\rho) = (p-aw)(v) = p(v)\ge 0.
$$
Conversely, let $q\in Q_y\cap \rho^\perp$. Then for any vertex
$v\in\Delta$ as before we have $q(v)= q(v-\rho)\ge 0$.
On the other hand, for any $v\in P_x\setminus S$ it holds $w(v)> 0$.
Therefore for $a\gg 0$ it holds $(q+aw)(v) \ge0$ for all $v\in P_x$.
Hence $q+aw\in P_x^\vee$ and $a=(q+aw)-aw\in (P_x^\vee)_{(w)}$.

Since $P_x^\vee$ is a sharp monoid, $0\in\Delta$ is a vertex. When we
view $Q_y^\vee$ as submonoid of $P_x^\gp$ as before, this vertex
corresponds to $-\rho\in Q_y^\vee$. It follows that $\RR_{\ge 0}\cdot
(-\rho) \subset Q_y^\vee$ is a one-dimensional face. Hence
$Q_y\cap \rho^\perp$ is a facet of $Q_y$, and $(\ref{identification
of cones})$ identifies $\Spec \big( k[(P_x^\vee)_{(w)}]\big)$ with
an irreducible component of the divisor $\chi^{\rho^*}=0$ in
$\Spec( k[Q_y])$:
$$
\Spec( k[(P_x^\vee)_{(w)}])\subset \Spec( k[Q_y]/(\chi^{\rho^*}).
$$
The left-hand side is a standard affine open subset of $\cY_x=
\Proj( k[P_x^\vee])$, which we denote by $U_{x,y}$. Hence
$$
\cY_x=\bigcup_y \cU_{x,y},
$$
where $y$ runs over all generic points $y\in F$ with $y\ge x$.
Moreover, by what we have just said, $\cU_{x,y}$ embeds as an
irreducible component into $\Spec k[Q_y]/(\chi^{\rho^*})$. 
Restriction also yields closed embeddings $\cU_{x,y}\cap \cY_z\to \Spec(
k[Q_y])$ for any $z\in F$ with $x\leq z\leq y$. These form a directed
system of closed embeddings, parametrized by all $z\leq y$. It is
compatible with the directed system defined by the $Y_x$ via a system
of open embeddings. Thus
\begin{eqnarray}\label{chart}
\cU_y:= \bigcup_{x\leq y} \cU_{x,y}\simeq \Spec( k[Q_y]/(\chi^{\rho^*})
\end{eqnarray}
is an open subscheme of $\cY_{(F,P,\rho)}$. The collection of the
closed embeddings $\cU_y\to \Spec k[Q_y]$, which is uniquely defined by
the duality datum, defines our pre-gtc atlas.

It remains to check the compatibility condition in the definition of
pre-gtc atlas (Definition~\ref{pre-gtc atlas} (ii)). Let $x\in F^*$ and
$q_x\in Q_x$. For generic points $y,y'\in F^*$ with $x\in \overline
y\cap \overline y'$ we have to show equality of the ideals $\shI$,
$\shI'$ on $\cU_y\cap \cU_{y'}$ generated by $q$ via the two gtc-charts
indexed by $y$ and $y'$. Denote $Q=Q_y$ and choose a lift $q\in Q$ of
$q_x$ under the generization map $Q\to Q_x$. It suffices to compare
the ideals on one of the open sets
$$
\cU_h=\Spec ( k[Q]_{(h)}),\quad h\in k[Q]
$$
generating the topology.

Let $v_1,\ldots,v_n$ be generators of the one-dimensional faces of
$Q^\vee$ and $\cU_i\subset \cU_h$ the irreducibe component corresponding
to $v_i$. By definition $I=(\chi^q)$. Precisely for those $i$ with $q(v_i)=0$
the function $\chi^q$ is non-zero at the generic point of $\cU_i$.
Therefore $\chi^q$ defines a Cartier divisor on the subspace
$\cZ\subset \cU_h$ corresponding to the ideal generated by
$$
\{p\in Q\,|\, q(v_i)=0\Rightarrow p(v_i)\neq 0 \text{ for all } i\}.
$$
The associated Weil divisor is
$
\sum_{q(v_j)\neq0} q(v_j)\cdot[\cZ\cap \cU_j]$.
The essential observation is that both $\cZ$ and this divisor depend
only on $q_x$. Hence, denoting by $f$ a generator of $I'$, there
exists $e\in k[Q]_{(h)}$, invertible on $\cZ$, with
$
f|_Z=(e\cdot \chi^q)|_\cZ$.
But $f$ and $\chi^q$ vanish at the generic points of the closure of
$\cU_h\setminus \cZ$, and hence $f=e\cdot \chi^q$ everywhere. This shows
$I'= I$.
\end{construction}

Next we describe the canonical involution on the set of all
duality data:

\begin{construction}\label{Mirror duality data}
(Mirror duality data.)
Let $(F,P,Q,\lambda)$ be a duality datum. The mirror duality datum
will be $(F^*,Q,P,\lambda^*)$, and we have to define the dual
compatibility datum $\lambda^*$. Let $\rho^*\in\Gamma(F^*,Q) $ be the
distinguished section, and let $x\leq y\in F$ be a closed and a generic
point, respectively. Recall that the given compatibility datum gives
an affine injection $\lambda_{xy}:Q_y^\vee\ra P^\gp_x$ with 
$\lambda_{xy}(0)=\rho$. The dual compatibility datum $\lambda^*$ is
defined by the formula
$$
\lambda^*_{yx}=((\lambda_{xy}-\rho)^\gp)^\vee +\rho^*:P_x^\vee\lra
Q_y^\gp.
$$
This indeed works:

\begin{lemma}
The collection $(F^*,Q,P,\lambda^*)$ is a duality datum.
\end{lemma}

\proof
We have to verify the compatibility condition
Definition~\ref{duality datum},(iii). Since $(\rho^*)^\perp\cap P_x$
is the facet belonging to $y$ we see that $-\rho^*$ spans the
one-dimensional face of $P_x^\vee$ corresponding to $y\in F^*$. Since
$\rho^\perp\cap Q_y$ is a facet of $Q_y$ and $P_x$ it remains to show
that
$$
\lambda^*_{yx}(P_x^\vee)\cap\rho^\perp\subset Q_y.
$$
For the following computation we view $P_x^\vee$ and $Q_y^\vee$ as
subsets of $Q_y^\gp$ and $P_x^\gp$ respectively. Let $m\in P_x^\vee$
with $(\rho^*+m)(\rho)=0$. We have to show that $(\rho^*+m)(p)\ge 0$
for any $p\in Q_y^\vee$. Since $Q_y^\vee$ is generated by elements of
the form $q-\rho$ with $q\in P_x\cap(\rho^*)^\perp$, we may restrict
to such elements. Now compute
$$
(\rho^*+m)(q-\rho)= \rho^*(q)-(\rho^*+m)(\rho)+m(q).
$$
The first two terms vanish, while $m(q)\ge 0$ since $m\in P_x^\vee$,
$q\in P_x$.
\qed
\end{construction}

It is clear from the definition of $\lambda^*$ that 
the mirror of the mirror $(F^*,Q,P,\lambda^*)$ is the original
duality datum $(F,P,Q,\lambda)$. In other words, passing to
the mirror duality datum defines an involution on the set of
duality data.

\section{Batyrev's mirror construction, degenerate abelian varieties}
In this section we illustrate our naive mirror construction with
two examples.

\begin{example}
(Batyrev's mirror construction.) Let $\Delta\subset \RR^n$ be a
polytope with integral vertices $v_i\in \ZZ^n$. We assume that
$\Delta$ is \emph{reflexive}, which means (1) the origin is the only
interior lattice point of $\Delta$, and (2) the polar polytope
$\Delta^\circ= \{m\in (\RR^n)^\vee\,|\, \langle m, v\rangle \ge -1\}$
has integral vertices. Then also the polar polytope is reflexive.
From $\Delta$ we obtain a duality datum as follows. Let $F$ be the
set of proper faces $\sigma\subsetneq\Delta$, where the relation
$\leq$ of points corresponds to inclusion  $\subset$ of faces. For
each face $\sigma$ define a monoid $P_\sigma$ as the quotient of the
``wedge monoid''
$$
\ZZ^n\cap \RR_{\ge0}\cdot\{p_2-p_1\,|\, p_1\in\sigma, p_2\in\Delta \}
$$
by its invertible elements. If $\sigma\subset\tau$ there is a
canonical surjection $P_\sigma\to P_\tau$ making these monoids into a
sheaf $P$ on $F$. Similarly, the polar polytope induces a sheaf $Q$
on the dual topological space $F^*$. For a face $\sigma\subset\Delta$
the monoid of integral points of the cone over $\sigma$ is
canonically dual to $Q_\sigma$. For every vertex $v\in\sigma$ we
therefore obtain an affine embedding $Q_\sigma^\vee \hookrightarrow
P_v$, and these provide the compatibility datum
Definition~\ref{duality datum} (iii). The Gorenstein property of both
$P$ and $Q$ follow from reflexivity of $\Delta$.

By going through the construction we see that $\cY_{(F,P,\rho)}$ is the
boundary divisor (the complement of the big cell) in the toric variety
$\PP(\Delta)$. The pre-gtc-atlas comes from the embedding into
$\PP(\Delta)$, so in this case actually glues to a logarithmic
structure. The conormal sheaf $\shN_\cD$ is the conormal sheaf of this
embedding. As it is never trivial, none of the global logarithmic structures
in the specified gtc-atlas is log-smooth over the standard log point.

The space $\cY_{(F^*,Q,\rho^*)}$ for the mirror duality datum gives the
boundary divisor in $\PP(\Delta^\circ)$. So here we retrieve part of
the Batyrev construction of mirror pairs of hypersurfaces in toric
varieties defined by reflexive polyhedra \cite{Batyrev 1994}. To go
further one would need to control the desingularization procedure
involved in Batyrev's construction under this process. This is beyond
the scope of this paper and will be further discussed in \cite{Gross;
Siebert 2003}.
\end{example}

\begin{example}
(Degenerate abelian varieties.)
Let $f:\ZZ^n\to \ZZ$ be a convex mapping, and $C_f\subset \RR^{n+1}$ 
the boundary of the convex hull of the graph
$
\Gamma_f=\big\{(v,f(v))\in\ZZ^{n+1}\,\big|\, v\in\ZZ^n\}.
$
Then $C_f$ is a multi-faceted paraboloid with integral vertices. We
assume all faces to be bounded. Let $F$ be the locally finite
topological space with points the faces of $C_f$ and the ordering
``$\leq$'' defined by inclusion of faces. Denote by
$\pi:\ZZ^{n+1}\to\ZZ^n$ the projection onto the $n$ first
coordinates. For each face $\sigma\subset C_f$ define
$$
Q_\sigma^\vee=\{(t \cdot v,t)\in \ZZ^{n+1}\,|\, v\in
\pi(\sigma),t\in\RR_{\ge 0}\}.
$$
As this ist the set of integral points of the cone over
$\pi(\sigma)$, embedded into the affine hyperplane $\{1\}\times
\RR^n$, there are compatible inclusions $Q_\sigma^\vee\lra
Q_\tau^\vee,\quad \sigma\subset\tau$. Therefore the duals $Q_\sigma$
form the stalks of a sheaf $Q$ on $F^*$. The projections
$Q_\sigma^\vee\to \NN$ onto the first coordinate define a section
$\rho^*$ of $Q$, and $(F^*,Q,\rho^*)$ is a gtc Kato fan.

Next we define the sheaf $P$ on $F$. By abuse of notation, for
$\sigma\in F$ let $\langle\sigma\rangle$ denote the saturated
subgroup of $\ZZ^{n+1}$ generated by $v-v'$ with $v,v'\in \sigma\cap
\ZZ^{n+1}$. Define
$$
P_\sigma\subset \ZZ^{n+1}/\langle \sigma\rangle
$$
to be the saturated submonoid generated by $w-v$, where $v\in\sigma
\cap \ZZ^{n+1}$ and $w\in C_f$, that is, $w= (\pi(w),t)$ with $t\ge
f(\pi(w))$. For $\sigma\subset\tau$ we have canonical surjections
$P_\sigma\ra P_\tau$, and this defines the sheaf $P$ on $F$. For
$\sigma\in F$ and $w\in\pi(\sigma)$ the equivalence class of
$(w,f(w)+1)$ in $P_\sigma$ defines the germ of the section $\rho$ at
$\sigma$. One can show that $(F^*,Q,\rho^*)$ is a gtc Kato fan.

For the compatibility datum let $v=(v_0,t_0)\in C_f$ be a
vertex and $\sigma\subset C_f$ a facet with $v\in\sigma$.
Then
$$
\lambda_{v\sigma}:Q_\sigma^\vee\lra P_v^\gp,\quad
(P,t)\longmapsto (P,t_0+1-t)
$$
is an affine embedding identifying $Q_\sigma^\vee$ with the integral
points of the cone over $\sigma$ with vertex
$\rho_\sigma=(v_0,t_0+1)$. $(F,P,Q,\Lambda)$ is a gtc-duality datum,
with $\cY_{(F^*,Q,\rho^*)}$ only \emph{locally} of finite type.

One can show that the mirror space $\cY_{(F,P,\rho)}$ is of the
same form, with defining function obtained by discrete Legendre
transform from $f$ \cite{Gross; Siebert 2002}, \cite{Gross; Siebert
2003}.

To obtain a degenerate abelian variety one assumes that $f=q+r$ with
$q(x)=x^tAx+b^tx+c$ a strictly convex quadratic function with
integral coefficients, and $r:\ZZ^n\to\ZZ$ a $\Lambda'$-periodic
function for a sublattice $\Lambda'\subset\Lambda:=\ZZ^n$ of finite
index. The $\Lambda'$-action on $\Lambda$ lifts to an affine action
on $\ZZ^{n+1}$ leaving $\Gamma_f$ invariant by
$$
T_w(v,\lambda)=(v+w,\lambda+2w^tAv+q(w)-c).
$$
The induced $\Lambda'$-action on the duality datum defines an
\'etale, quasicompact equivalence relation on $\cY_{(F^*,Q,\rho^*)}$.
The quotient $\cY_{(F^*,Q,\rho^*)}/\Lambda'$ is the central fiber of
the degeneration of polarized abelian varieties associated to $q+r$
by Mumford's construction \cite{Mumford 1972}. The quotient of the
gtc-atlas gives the log structure associated to the degeneration. So
here there actually is a log-smooth morphism to the standard log
point. Up to changing the gluing of the irreducible components any
maximally degenerate polarized abelian variety is of this form
\cite{Alexeev 2002}, Section~5.7. In the mirror picture $\Lambda^*$
is the sublattice of $\Lambda^\vee$ generated by the slopes of $f$,
while $(\Lambda^*)'=\Lambda'$ with action induced from the action on
$\Gamma_f$.

For an explicit two-dimensional example take $\Lambda'=2\ZZ^2$,
$q(x,y)= x^2-xy+y^2$ and $r(v)=1$ for $v\in\Lambda'$ and $r(v)=0$
otherwise. Then $\cY_{(F^*,Q,\rho^*)}/\Lambda'$ is a union of $3$ copies of
$\PP^1\times \PP^1$ and in each copy, the pull-back of the singular
locus is a $4$-gon of lines. The mirror $\cY_{(F,P,\rho)}/
(\Lambda^*)'$ is a union of $2$ copies of $\PP^2$ and a $\PP^2$ blown
up in $3$ points. The pull-back of the singular locus is a union of
$3$ lines for $\PP^2$, and a $6$-gon of rational curves containing
the exceptional curves for the other component.
\end{example}


\end{document}